\documentclass[12pt,a4paper]{article}

\usepackage{amssymb}

\usepackage{amsmath}

\usepackage{cases}

\newtheorem{theorem}{Theorem}[section]
\newtheorem{definition}[theorem]{Definition}
\newtheorem{lemma}[theorem]{Lemma}
\newtheorem{corollary}[theorem]{Corollary}
\newtheorem{proposition}[theorem]{Proposition}
\newtheorem{example}[theorem]{Example}
\newtheorem{remark}[theorem]{Remark}

\usepackage{cite}

\begin{document}
\title{A new Composition-Diamond lemma for dialgebras\footnote{Supported by the NNSF of China (11571121) and the Science and Technology Program of Guangzhou (201605121833161).}}
\author{Guangliang Zhang and
 Yuqun
Chen\footnote {Corresponding author.}  \\
{\small \ School of Mathematical Sciences, South China Normal
University}\\
{\small Guangzhou 510631, P. R. China}\\
{\small  yqchen@scnu.edu.cn}\\
{\small zgl541@163.com}}

\date{}

\maketitle \noindent\textbf{Abstract:} Let $Di\langle X\rangle$ be the free dialgebra over a field generated by a set $X$.
Let $S$ be a monic subset of $Di\langle X\rangle$.  A Composition-Diamond lemma for dialgebras is
firstly established by Bokut, Chen and Liu in 2010 \cite{Di}
which claims that if (i) $S$ is a Gr\"{o}bner-Shirshov basis in $Di\langle X\rangle$, then (ii) the
set of $S$-irreducible words is a linear basis of the quotient dialgebra $Di\langle X \mid S \rangle$,
but not conversely.
Such a lemma based on  a fixed ordering on normal diwords  of $Di\langle X\rangle$ and special definition of composition trivial modulo $S$.
In this paper, by
introducing an arbitrary monomial-center ordering and the usual definition of composition  trivial modulo $S$, we give a new Composition-Diamond lemma for dialgebras which makes the conditions (i) and (ii) equivalent. We show that every ideal of $Di\langle X\rangle$ has a unique reduced Gr\"{o}bner-Shirshov basis.
The new lemma is more useful and convenient than the one in \cite{Di}.
As applications, we give a method to find  normal forms of elements of an arbitrary disemigroup, in particular, A.V. Zhuchok's (2010) and Y.V. Zhuchok's (2015) normal forms of the free commutative disemigroups and the free abelian disemigroups, and normal forms of the free left (right) commutative disemigroups.



\section{Introduction}

The notion of a dialgebra (disemigroup) was introduced by Loday \cite{Lo95} and
investigated in many papers (see, for example, \cite{Di,Fe,Ko,Lo95,Po,Zhuchok,Zhuchok17,Zhu15}).
Loday \cite{Lo95} constructed a free dialgebra and the universal enveloping dialgebra for a Leibniz algebra.
Bokut, Chen and Liu \cite{Di} established Gr\"{o}bner-Shirshov bases theory for dialgebras.
Pozhidaev \cite{Po} studied the connection of Rota-Baxter algebras and dialgebras with associative bar-unity.
Kolesnikov \cite{Ko} proved recently that each dialgebra may be obtained in turn from an associative conformal algebra.
Analogues of some notions of the functional analysis were defined on dialgebras in \cite{Fe}. A.V. Zhuchok \cite{Zhuchok} and Y.V. Zhuchok \cite{Zhu15} constructed the free commutative disemigroup and the free abelian disemigroup respectively. Various free disemigroups were introduced by A.V. Zhuchok in a survey paper \cite{Zhuchok17}.

Gr\"{o}bner bases and Gr\"{o}bner-Shirshov bases were invented independently by A.I. Shirshov
for ideals of free (commutative, anti-commutative) non-associative algebras \cite{Sh62b,Shir3},
free Lie algebras \cite{Shir3} and implicitly free associative algebras \cite{Shir3} (see also \cite{be78,bo76}), by
H. Hironaka \cite{Hi64} for ideals of the power series algebras (both formal and convergent),
and by B. Buchberger \cite{bu70} for ideals of the polynomial algebras. Gr\"{o}bner bases and
Gr\"{o}bner-Shirshov bases theories have been proved to be very useful in different branches
of mathematics, including commutative algebra and combinatorial algebra. It is a powerful
tool to solve the following classical problems: normal form; word problem; conjugacy
problem; rewriting system; automaton; embedding theorem; PBW theorem; extension;
homology; growth function; Dehn function; complexity; etc. See, for example, the books
\cite{AL,BKu94,BuCL,BuW,CLO,Ei} and the surveys \cite{BC,BokutChenBook,BC14,BoFKK00,BK03,BK05a}.

In Gr\"{o}bner-Shirshov bases theory for a category of algebras, a key part is to establish
``Composition-Diamond lemma" for such algebras. The name ``Composition-Diamond
lemma" combines the Neuman Diamond Lemma \cite{Nm}, the Shirshov Composition Lemma \cite{Sh62b}
 and the Bergman Diamond Lemma \cite{be78}. A Composition-Diamond lemma for
dialgebras was firstly given
by Bokut, Chen and Liu in 2010 \cite{Di}.

Let $Di\langle X\rangle$ be the free dialgebra over a field $\mathbf{k}$ generated by a well-ordered set $X$ and $X^+$ the free semigroup generated by $X$ without the unit. With the notation as in \cite{Di}, for any $u=x_1\cdots x_{m}\cdots x_n\in X^+$,
$$
[u]_m:=x_1\cdots x_{m-1} \dot{x_m} x_{m+1}\cdots x_n
=x_1\vdash \cdots\vdash x_{m-1}\vdash x_m\dashv x_{m+1}\dashv \cdots \dashv x_n 
$$
is called a  normal diword on $X$. The set $[X^+]_\omega$ of all  normal diwords on $X$ is a linear basis of $Di\langle X\rangle$.
Let $[X^+]_\omega$ be a well-ordered set, $S\subset Di\langle X\rangle$  a monic subset
of polynomials and $Id(S)$ be the ideal of $Di\langle X\rangle$ generated by $S$. A normal diword $[u]_n$ is
said to be $S$-irreducible if $[u]_n$ is not equal to the leading monomial of any normal $S$-diword. Let Irr$(S)$ be
the set of all $S$-irreducible diwords. Consider the following statements:

(i) The set $S$ is a Gr\"{o}bner-Shirshov basis in $Di\langle X\rangle$.

(ii) The set Irr$(S)$ is a $\mathbf{k}$-basis of the quotient dialgebra  $Di\langle X \mid S\rangle:= Di\langle X\rangle/Id(S)$.

In \cite{Di}, it is shown that $(i) \Rightarrow(ii)$ but $(ii)\nRightarrow (i)$.
Their proof of the above result based on  a fixed ordering on $[X^+]_\omega$ and special definition of composition trivial modulo $S$.
In this paper, for an arbitrary monomial ordering on $X^+$, we introduce a so-called monomial-center ordering on $[X^+]_\omega$ and give a new Composition-Diamond lemma for dialgebras
which makes the two conditions above equivalent, see Theorem \ref{cd}. Comparing with the corresponding result in \cite{Di}, the new lemma will be more useful and convenient  when one calculates a Gr\"{o}bner-Shirshov basis in $Di\langle X\rangle$. We show that with a monomial-center ordering, every ideal of $Di\langle X\rangle$ has a unique reduced Gr\"{o}bner-Shirshov basis.
As applications, we give a method to find normal forms of elements of an arbitrary disemigroup.
In particular, we give short proofs of A.V. Zhuchok \cite{Zhuchok} and Y.V. Zhuchok's \cite{Zhu15} results on normal forms of elements of  the free commutative disemigroup and the free abelian disemigroup generated by a set $X$, respectively.
Moreover, we give Gr\"{o}bner-Shirshov bases for some dialgebras and disemigroups,
and obtain  normal forms of elements of them.

The paper is organized as follows. In section 2, we review the free dialgebra $Di\langle X\rangle$ over a field $\mathbf{k}$ generated by $X$.
In section 3, by introducing a monomial-center ordering on $[X^+]_\omega$, normal $S$-diwords and compositions, we give a new Composition-Diamond lemma for dialgebras which makes the conditions (i) and (ii) mentioned before equivalent. In section 4, Gr\"{o}bner-Shirshov bases theory for dirings is introduced, which may find an $R$-basis for some disemigroup-dirings over an associative ring $R$. In section 5,  some applications are given.

\ \

\section{Preliminaries}

Throughout the paper, we fix a field $\mathbf{k}$.
$\mathbb{Z}^+$ stands for the set of positive integers.
For a nonempty set $X$, we define the following notations:

$X^*$: the set of all associative words on $X$ including the empty word, i.e. the free monoid generated by $X$.

$X^+$: the set of all nonempty associative words on $X$, i.e. the free semigroup generated by $X$ without the unit.

$\lfloor X^+\rfloor:= \{\lfloor x_{i_1}x_{i_2}\cdots  x_{i_n}\rfloor\mid  i_1,\dots,i_n\in I, i_1\leq i_2\leq \cdots \leq i_n, n\in \mathbb{Z}^+\}$, the set of all nonempty commutative associative words on $X$, where $X=\{x_i\mid i\in I\}$ is a total-ordered set.

$[X^+]_\omega:=\{[u]_m\mid u\in X^+, m\in \mathbb{Z}^+, 1\leq m \leq |u|\}$, the set of all associative normal diwords on $X$, following the notation in \cite{Di}, where $|u|$ is the number of letters in $u$ (the length of $u$).

$\lfloor X^+\rfloor_\omega:=\{\lfloor u\rfloor_m\mid \lfloor u\rfloor\in \lfloor X^+\rfloor, m\in \mathbb{Z}^+, 1\leq m \leq |u|\}$, the set of all commutative normal diwords on $X$.

\ \ For $ u\in X^+,\ [u]_m$ is called an associative diword, while $\lfloor u\rfloor_m$ is called a commutative diword. For example, if $u=x_2x_1x_2x_1\in X^+,\ x_1<x_2$, then $\lfloor u\rfloor=\lfloor x_1x_1 x_2x_2\rfloor, \ [u]_3=x_2x_1 \dot{x_2}x_1,\  \lfloor u\rfloor_3=\lfloor x_1x_1x_2x_2\rfloor_3=x_1x_1 \dot{x_2}x_2$.

$\lfloor X^+ \rfloor_{_1}:=\{\lfloor u \rfloor_1 \mid \lfloor u\rfloor\in \lfloor X^+\rfloor\}$.

$\lfloor X^+ \rfloor_{_{2-2}}:=\{\lfloor v \rfloor_2 \mid \lfloor v\rfloor\in \lfloor X^+\rfloor, |v|=2\}$.

$Di\langle X\rangle$: the free dialgebra over a field $\mathbf{k}$ generated by $X$.

$Di_R\langle X\rangle$: the free diring over an associative ring $R$ generated by $X$.

$Disgp\langle X\rangle=[X^+]_\omega$: the free disemigroup generated by $X$.

$Di[ X]$: the free commutative dialgebra over a field $\mathbf{k}$ generated by $X$.

$Disgp[ X]=\lfloor X^+ \rfloor_{_1}\cup\lfloor X^+ \rfloor_{_{2-2}}$: the free commutative disemigroup generated by $X$.

\begin{definition}(\cite{Lo99})\label{l2}
\emph{An} associative dialgebra \emph{(}dialgebra \emph{for short) is a} $\mathbf{k}$-\emph{module} $D$ \emph{equipped with two}
$\mathbf{k}$-\emph{linear maps}
$$
\vdash \ : D\otimes D\rightarrow D,
\ \ \ \
\dashv \ : D\otimes D\rightarrow D,
$$
\emph{where} $\vdash$ \emph{and} $\dashv$ \emph{are associative and satisfy the following identities:}
\begin{equation}\label{eq00}
\begin{cases}
a\dashv(b\vdash c)=a\dashv (b\dashv c),\\
(a\dashv b)\vdash c=(a\vdash b)\vdash c, \\
a\vdash(b\dashv c)=(a\vdash b)\dashv c.
\end{cases}
\end{equation}
\emph{for all} $a, \ b, \ c\in D$.

\emph{A dialgebra} $(D,\vdash,\dashv)$ \emph{is}  commutative \emph{if both } $\vdash$ \emph{and} $\dashv$ \emph{are commutative}.

\end{definition}

Write
$$
[X^+]_\omega:=\{[u]_m \mid u\in X^+, m\in \mathbb{Z}^+, 1\leq m \leq |u|\},
$$
where $|u|$ is the number of letters in $u$. For any $h=[u]_m \in [X^+]_\omega$,
we call $u$ the \textit{associative word} of $h$, and
$m$, denoted by $p(h)$, the \textit{position of center} of $h$.
For example, if $u=x_1x_2\cdots x_n\in X^+$, $x_i\in X,\ h=[u]_m,\ 1\leq m\leq n$, then $p(h)=m$ and with the notation as in \cite{Di},
$$
[u]_m:=x_1\cdots x_{m-1} \dot{x_m} x_{m+1}\cdots x_n
=x_1\vdash \cdots\vdash x_{m-1}\vdash x_m\dashv x_{m+1}\dashv \cdots \dashv x_n.
$$
A  word $[u]_m\in [X^+]_\omega$ is called a \textit{normal diword}.

Let $Di\langle X\rangle$ be the free $\mathbf{k}$-module with a $\mathbf{k}$-basis $[X^+]_\omega$.
For any $[u]_m, [v]_n \in [X^+]_\omega$, define
$$
[u]_m\vdash [v]_n=[uv]_{|u|+n}, \ \ \ [u]_m\dashv [v]_n=[uv]_m,
$$
and extend them linearly to $Di\langle X\rangle$.
It is well known from \cite{Lo99} that $Di\langle X\rangle$ is the free dialgebra generated by $X$.

Let $X$ be a well-ordered set. We define the \textit{deg-lex ordering} on $X^+$ by the following:
for $u=x_{i_1}x_{i_2}\cdots x_{i_n}, v=x_{j_1}x_{j_2}\cdots x_{j_m}\in X^+$, where each $x_{i_l},x_{j_t}\in X$,
\begin{equation*}\label{equ0}
u>v \ \Leftrightarrow \ (|u|,x_{i_1},x_{i_2},\cdots, x_{i_n})>(|v|,x_{j_1},x_{j_2},\cdots, x_{j_m}) \ \mbox{lexicographically}.
\end{equation*}
An ordering $>$ on $X^+$ is said to be \textit{monomial} if $>$ is a well ordering and for any $u,v,w\in X^+$,
$$
u>v \Rightarrow uw>vw\ \mbox{ and }\ wu>wv.
$$
Clearly, the deg-lex ordering is monomial.

\section{A new Composition-Diamond lemma}

Let $>$ be a monomial ordering on $X^+$.
We define the \textit{monomial-center ordering} $>_d$ on $[X^+]_\omega$ as follows. For any $[u]_m,[v]_n\in [X^+]_\omega$,
\begin{equation}\label{equ0}
[u]_m>_d[v]_n \ \mbox{if} \  (u,m)>(v,n) \ \ \mbox{lexicographically}.
\end{equation}
In particular, if $>$ is the deg-lex ordering on $X^+$, we call the ordering defined by $(\ref{equ0})$ the \textit{deg-lex-center ordering} on $[X^+]_\omega$.
For simplicity of notation, we write $>$ instead of $>_d$ when no confusion can arise.
It is clear that a monomial-center ordering is a well ordering on $[X^+]_\omega$. Such an ordering plays an important role in this paper. Here and subsequently, the monomial-center ordering on $[X^+]_\omega$ will be used, unless otherwise stated.

For convenience we assume that $[u]_m>0$ for any $[u]_m\in [X^+]_\omega$.
For any nonzero polynomial $f\in Di\langle X\rangle$, let us denote $\overline{f}$ be the \textit{leading monomial} of $f$ with respect to the ordering $>$,
$lt(f)$ the \textit{leading term} of $f$, $lc(f)$ the \textit{coefficient} of $\overline{f}$ and $\widetilde{f}$ the associative word of $\overline{f}$.
$f$ is called \textit{monic} if $lc(f)=1$. For any nonempty subset $S$ of $Di\langle X\rangle$,
$S$ is\textit{ monic} if $s$ is monic for all $s\in S$.

\begin{definition}
\emph{A nonzero polynomial} $f\in Di\langle X\rangle$ \emph{is} strong \emph{if} $\widetilde{f}> \widetilde{r\!_{_f}}$, \emph{where} $r\!_{_f}:=f-lt(f)$.
\end{definition}

It is easy to check that
$>$ on $[X^+]_\omega$ is monomial in the following sense:
\begin{eqnarray*}{ }
\ [u]_m>[v]_n &\Rightarrow &[w]_l \vdash [u]_m  > [w]_l \vdash [v]_n,\\
&& [u]_m \dashv [w]_l >  [v]_n \dashv [w]_l,\\
&&  [u]_m \vdash [w]_l \geq [v]_n \vdash [w]_l, \\
&& [w]_l \dashv [u]_m \geq [w]_l \dashv [v]_n, \\
\ u>v  &\Rightarrow& [u]_m \vdash [w]_l > [v]_n \vdash [w]_l, \\
 &&[w]_l \dashv [u]_m > [w]_l \dashv [v]_n,
 \end{eqnarray*}
where $[u]_m, [v]_n,  [w]_l\in [X^+]_\omega$.

From this it follows that
\begin{lemma}\label{r1}
Let $0\neq f\in Di\langle X\rangle$ and $[u]_m \in [X^+]_\omega$. Then
 \begin{eqnarray*}
\overline{([u]_m\vdash f)}=[u]_m \vdash \overline{f},&& \ \  \ \overline{(f\dashv [u]_m)}= \overline{f}\dashv [u]_m,\\
 \overline{([u]_m\dashv f)} \leq [u]_m \dashv \overline{f},&& \ \ \ \overline{(f\vdash [u]_m)}\leq \overline{f}\vdash [u]_m.
 \end{eqnarray*}
In particular, if $f$ is strong, then $\overline{([u]_m\dashv f)} = [u]_m \dashv \overline{f}$ and $\overline{(f\vdash [u]_m)}= \overline{f}\vdash [u]_m$.
\end{lemma}

\begin{example}
Let $X=\{x_1,x_2,x_3\}$, $x_1>x_2>x_3$, $Char\mathbf{k}\neq2,3$ and $>$ be the deg-lex-center ordering on $[X^+]_\omega$.
Let $f=2[x_1x_2x_3]_3-2[x_1x_2x_3]_2+3[x_1x_3]_2$. Then
$$
\overline{f}=[x_1x_2x_3]_3,\ lt(f)=2[x_1x_2x_3]_3,\ lc(f)=2,\ \widetilde{f}=x_1x_2x_3, \ r\!_{_f}=-2[x_1x_2x_3]_2+3[x_1x_3]_2.
$$
The polynomial $f$ is not strong since $\widetilde{f}=x_1x_2x_3=\widetilde{r\!_{_f}}$. Of course, $r\!_{_f}$ is strong. It is easy to check that
 \begin{eqnarray*}
\overline{(x_1\vdash f)}&=&[x_1x_1x_2x_3]_4=x_1\vdash \overline{f}, \\ \overline{(f \dashv x_1)}&=&[x_1x_2x_3x_1]_3=\overline{f}\dashv x_1,\\
\overline{(x_1\dashv f)}&=&[x_1x_1x_3]_1<[x_1x_1x_2x_3]_1=x_1\dashv \overline{f}, \\
\overline{(f \vdash x_1)}&=&[x_1x_3x_1]_3<[x_1x_2x_3x_1]_4=\overline{f}\vdash x_1,\\
\overline{(x_1\dashv r\!_{_f})}&=&[x_1x_1x_2x_3]_1=x_1\dashv \overline{r\!_{_f}},\\
\overline{(r\!_{_f} \vdash x_1)}&=&[x_1x_2x_3x_1]_4=\overline{r\!_{_f}}\vdash x_1.
\end{eqnarray*}
\end{example}

Here and subsequently, $S$ denotes a monic subset of $Di\langle X\rangle$ unless otherwise stated.

By an $S$-\textit{diword} $g$ we mean
a normal diword on $X\cup S$ with only one occurrence of $s \in S$. If this is the case and
\begin{eqnarray}\label{d1}
g=[x_{i_1}\cdots x_{i_k}\cdots x_{i_n}]_m|_{_{x_{i_k}\mapsto s}},
\end{eqnarray}
where $1\leq k \leq n,\ x_{i_l}\in X,\ 1\leq l \leq n$, then we also call $g$ an \textit{$s$-diword}. For simplicity, we denote the $s$-diword of the form (\ref{d1}) by $(asb)$, where $a,b\in X^*$, $s\in S$.

\begin{definition}\label{dns}
\emph{An} $S$-\emph{diword} \emph{(\ref{d1})} \emph{is called a} normal $S$-diword \emph{if either} $k=m$ \emph{or} $s$ \emph{is strong}.
\end{definition}

Note that if $(asb)$ is a normal $S$-diword, then $\overline{(asb)}=[a\widetilde{s}b]_l$ for some $l\in P([asb])$, where
\begin{displaymath} P([asb]):=
\begin{cases}
\{n\in \mathbb{Z}^+\mid 1\leq n \leq |a|\} \cup \{|a|+p(\overline{s})\} \cup \{n\in \mathbb{Z}^+\mid |a\widetilde{s}|<n\leq |a\widetilde{s}b|\}
 \ \ \text{if $s$ is strong,}\\
\{|a|+p(\overline{s})\}  \ \  \text{if $s$ is not strong.}
\end{cases}
\end{displaymath}
If this is so, we denote the normal $S$-diword $(asb)$ by $[asb]_l$ and  also call $[asb]_l$ a \textit{normal $s$-diword}.

In what follows, to simplify notation, we let
\begin{displaymath} [u]_m\vdash f \dashv [v]_n:=
\begin{cases}
[u]_m\vdash f
 \ \ \text{if $v$ is empty,}\\
f \dashv [v]_n  \ \  \text{if $u$ is empty,}
\end{cases}
\end{displaymath}
where $[u]_m,[v]_n\in [X^+]_\omega$, $f\in Di\langle X\rangle$.
The lemma below follows from Definition \ref{dns}.

\begin{lemma}\label{l0}
Let $(asb)$ be an $s$-diword and $[u]_m,[v]_n\in [X^+]_\omega$.
Then $(asb)=[asb]_l$
if and only if $[u]_m\vdash (asb) \dashv [v]_n=[uasbv]_{|u|+l}$, where $u,v$ may be empty.
\end{lemma}

\begin{definition}\label{dcom}
\emph{Let} $f,g$ \emph{be monic polynomials in} $Di\langle X\rangle$.
\begin{enumerate}
\item[1)] \emph{If} $f$ \emph{is not strong, then we call} $x \dashv f$ \emph{the} composition of left multiplication \emph{of} $f$
\emph{for all} $x\in X$ \emph{and} $f \vdash [u]_{|u|}$ \emph{the} composition of right multiplication \emph{of} $f$  \emph{for all} $u\in X^+$.

\item[2)] \emph{Suppose that}  $w=\widetilde{f}= a\widetilde{g}b$ \emph{for some} $a,b\in X^*$
\emph{and} $(agb)$ \emph{is a normal} $g$-\emph{diword}.

2.1 \ \emph{If} $p(\overline{f})\in P([agb])$,
\emph{then we call}
$$
(f,g)_{\overline{f}}=f-[agb]_{p(\overline{f})}
$$
\emph{the} composition of inclusion \emph{of} $f$ \emph{and} $g$.

2.2 \ \emph{If} $p(\overline{f})\notin P([agb])$ \emph{and both} $f$ \emph{and} $g$ \emph{are strong},
\emph{then for any} $x\in X$ \emph{we call}
$$
(f,g)_{[xw]_1}=[xf]_1-[xagb]_1
$$
\emph{the} composition of left multiplicative inclusion \emph{of} $f$ \emph{and} $g$, \emph{and}
$$
(f,g)_{[wx]_{_{|wx|}}}=[fx]_{|wx|}-[agbx]_{|wx|}
$$
\emph{the} composition of right multiplicative inclusion \emph{of} $f$ \emph{and} $g$.

\item[3)] \emph{Suppose that there exists a} $w=\widetilde{f}b= a\widetilde{g}$ \emph{for some} $a,b\in X^*$
\emph{such that} $|\widetilde{f}|+|\widetilde{g}|>|w|$,  $(fb)$ \emph{is a normal $f$-diword
and} $(ag)$ \emph{is a normal} $g$-\emph{diword}.

3.1 \ \emph{If} $P([fb])\cap P([ag])\neq \varnothing$,
\emph{then for any} $m\in P([fb])\cap P([ag])$ \emph{we call}
$$
(f,g)_{[w]_m}=[fb]_m-[ag]_m
$$
\emph{the} composition of intersection \emph{of} $f$ \emph{and} $g$.

3.2 \ \emph{If} $P([fb])\cap P([ag])=\varnothing$ \emph{and both} $f$ \emph{and} $g$ \emph{are strong},
\emph{then for any} $x\in X$ \emph{we call}
$$
(f,g)_{[xw]_1}=[xfb]_1-[xag]_1
$$
\emph{the} composition of left multiplicative intersection \emph{of} $f$ \emph{and} $g$, \emph{and}
$$
(f,g)_{{[wx]}_{|wx|}}=[fbx]_{|wx|}-[agx]_{|wx|}
$$
\emph{the} composition of right multiplicative intersection \emph{of} $f$ \emph{and} $g$.
\end{enumerate}
\end{definition}

For any composition $(f,g)_{[u]_n}$ mentioned above, we call $[u]_n$ the \textit{ambiguity} of $f$ and $g$.

\begin{definition}\label{dgsb}
\emph{Let} $S$ \emph{be  a monic subset of} $Di\langle X\rangle$ \emph{and} $[w]_m\in [X^+]_\omega$.
\emph{A polynomial} $h\in Di\langle X\rangle$ \emph{is} trivial modulo $S$ \emph{(}$(S,[w]_m)$, \emph{resp.)}, \emph{denoted by}
$$
h\equiv 0 \ mod(S) \  (mod(S,[w]_m), \emph{resp}.),
$$
\emph{if} $h=\sum{\alpha}_i[a_i s_i b_i]_{m_i}$, \emph{where each} $\alpha_i\in \mathbf{k}, \ a_i,b_i\in X^*, \ s_i\in S$
\emph{and} $\overline{[a_i s_ib_i]_{m_i}}\leq \overline{h}$ \emph{(}$\overline{[a_i s_ib_i]_{m_i}}<[w]_m$, \emph{resp}.\emph{)}.

\emph{A monic set} $S$ \emph{is called a} Gr\"{o}bner-Shirshov basis \emph{in} $Di\langle X\rangle$ \emph{if any
composition of polynomials in} $S$ \emph{is trivial modulo} $S$.
\end{definition}

A monic set $S$ is said to be \textit{closed} under the composition of left (right, resp.) multiplication
if all left (right, resp.) multiplication compositions of elements of $S$ are trivial modulo $S$.
We set
\begin{enumerate}
\item[] \
$Irr(S):=\{[u]_n\in [X^+]_\omega\mid [u]_n\neq \overline{[asb]_m}$ for any normal $S$-diword $[asb]_m$ \}.
\end{enumerate}

\begin{remark}\label{rcomp}
\emph{The definition of a Gr\"{o}bner-Shirshov basis in \cite{Di}
is different from the Definition \ref{dgsb}.} \emph{In \cite{Di}, the definition of a Gr\"{o}bner-Shirshov basis is based on a fixed ordering on $[X^+]_\omega$.}
\emph{Comparing with \cite{Di}, we have different definitions of the following: ordering of normal diwords; normal $S$-diword;
compositions of left and right multiplication,  multiplicative inclusion and multiplicative intersection; and composition to be trivial.}
\end{remark}

In calculating a Gr\"{o}bner-Shirshov basis in $Di\langle X\rangle$, the following example shows that our method is more convenient than the one of \cite{Di}. In all examples of this section, we let $>$ be the deg-lex-center ordering on $[X^+]_\omega$, where $X$ is a well-ordered set.

\begin{example}
Let $D=Di\langle X\mid S\rangle$.
If $S\subseteq [X^+]_\omega$, then it is easy to check that $S$ is a Gr\"{o}bner-Shirshov basis.
But the result is not true in the sense of \cite{Di}.
For example, let $X=\{x_1, x_2, x_3\}$, $x_1>x_2>x_3$, and $D=Di\langle X\mid [x_1x_2]_2\rangle$.
Then $S=\{[x_1x_2]_2\}$ is a Gr\"{o}bner-Shirshov basis.
Applying Theorem \ref{cd} we conclude that
$$
Irr(S)=\{[z_m\dots z_1xy_1\dots y_n]_{m+1}\mid z_j,x,y_i\in X,
z_{j+1}z_j\neq x_1x_2, y_iy_{i+1}\neq x_1x_2, z_1x\neq x_1x_2\}
$$
is a linear basis of $D$. Let $S_1=\{[x_1x_2]_2,\ [xx_1x_2]_1 \mid x\in X\}$.
In the sense of \cite{Di}, $S_1$ is a Gr\"{o}bner-Shirshov basis, but $S$ is not.
However, $Irr(S_1)$ in the sense of \cite{Di} is the same as the set $Irr(S)$.
\end{example}

\begin{lemma}\label{lnew}
Let $S$ be closed under the composition of left multiplication and $f\in S$.
If $f$ is not strong, then for any $[u]_1 \in [X^+]_\omega$,
$[u]_1\dashv f \equiv 0 \ mod(S)$.
\end{lemma}
\noindent{\bf Proof.} The proof follows by induction on $(u\widetilde{f},|u|)$.
If $|u|=1$, then the result holds.
Assume that $|u|\geq 2$ and $[u]_1=[vx]_1$, $v\in X^+$, $x\in X$. Then
$[u]_1\dashv f=[v]_1\dashv (x\dashv f)$ is a linear combination of $S$-diwords of the form
$
[v]_1\dashv [asb]_m,
$
where $s\in S$ and $[a\widetilde{s}b]_m\leq \overline{(x\dashv f)}$.
It follows that $\overline{([v]_1\dashv [asb]_m)}\leq [v]_1\dashv [a\widetilde{s}b]_m \leq [v]_1\dashv \overline{(x\dashv f)}=\overline{([u]_1\dashv f)}$ and $a\widetilde{s}b\leq x\widetilde{f}$.

If $s$ is strong, then $[v]_1\dashv [asb]_m$ is already a normal $S$-diword, and we have done.

Suppose that $s$ is not strong. If $a$ is empty, then $[v]_1\dashv [asb]_m=([v]_1\dashv s)\dashv [b]_1$ and $(v\widetilde{s},|v|)<(u\widetilde{f},|u|)$.
If $a$ is not empty, then $[v]_1\dashv [asb]_m=([va]_1\dashv s)\dashv [b]_1$ and $m=|a|+p(\overline{s})>1$.
Since $[a\widetilde{s}b]_m\leq [x\widetilde{f}]_1$, we have $a\widetilde{s}b<x\widetilde{f}$ and $(va\widetilde{s},|va|)<(u\widetilde{f},|u|)$.
By induction, $[v]_1\dashv [asb]_m$ is a linear combination of $S$-diwords of the form
$[cs'd]_n\dashv [b]_1$, where $s'\in S$ and $[c\widetilde{s}'d]_n\leq \overline{([va]_1\dashv s)}$.
By Lemma \ref{l0}, $[cs'd]_n\dashv [b]_1$ is a normal $S$-diword, and
$[c\widetilde{s}'d]_n\dashv [b]_1\leq \overline{([va]_1\dashv s)}\dashv [b]_1=\overline{([v]_1\dashv [asb]_m)}\leq \overline{([u]_1\dashv f)}$.
\ \ $\square$

\begin{remark}
\emph{The following example shows that
Lemma \ref{lnew}  is not true if we replace} ``$x\dashv f$'',  ``$[u]_1\dashv f$'' \emph{by} ``$f \vdash x$'', ``$f \vdash [u]_{|u|}$'' \emph{respectively}.
\end{remark}

\begin{example}
Let $X=\{x_1, x_2\}, x_1>x_2$, $Char\mathbf{k}\neq2$ and $S=\{f, g, h \}$, where $f=[x_1x_2]_2+[x_1x_2]_1$,
$g=[x_1x_2x_1]_3-\frac{1}{2}[x_1x_2x_1]_2-\frac{1}{2}[x_1x_2x_1]_1$,
$h=[x_1x_2x_2]_3-\frac{1}{2}[x_1x_2x_2]_2-\frac{1}{2}[x_1x_2x_2]_1$.
Clearly, $f,g$ and $h$ are not strong. We check at once that $g\vdash x_i=0$, $h\vdash x_i=0$, $i=1,2$, and
$$
f\vdash x_1=2g+f\dashv x_1\equiv 0 \ mod(S), \ \ f\vdash x_2=2h+f\dashv x_2\equiv 0 \ mod(S).
$$
However, $f\vdash [x_1x_1]_2$ is not trivial modulo $S$.
\end{example}

\begin{lemma}\label{lnewa}
Let $S$ be closed under the compositions of left and right multiplication.
Then for any normal $S$-diword $[asb]_m$ and $[u]_n \in [X^+]_\omega$,
$$
[u]_n\dashv [asb]_m\equiv 0 \ mod(S), \ \ [asb]_m\vdash [u]_n \equiv 0 \ mod(S).
$$
Moreover, if $a\widetilde{s}b<w$, $w\in X^+$,
then
$$
[u]_n\dashv [asb]_m\equiv 0 \ mod(S,\ [uw]_n), \ \ [asb]_m\vdash [u]_n \equiv 0 \ mod(S, \ [wu]_{|w|+n}).
$$
\end{lemma}
\noindent{\bf Proof.}
We prove only the results for the case $[u]_n\dashv [asb]_m$. The proof of the another case is similar.

If $s$ is strong, then $[u]_n\dashv [asb]_m$ is a normal $S$-diword.
Assume that $s$ is not strong. Note that $[u]_n=[u_1]_{|u_1|}\vdash[u_2]_1$, $u_1,u_2\in X^*, |u_1|=n-1$.
Then
$$
[u]_n\dashv [asb]_m=[u_1]_{|u_1|}\vdash ([u_2a]_1\dashv s)\dashv [b]_1.
$$
By Lemma \ref{lnew},
$[u]_n\dashv [asb]_m$ is a linear combination of $S$-diwords of the form
$
[u_1]_{|u_1|}\vdash[cs'd]_l\dashv [b]_1,
$
where $s'\in S$ and $[c\widetilde{s}'d]_l\leq \overline{([u_2a]_1\dashv s)}$.
By Lemma \ref{l0}, $[u_1]_{|u_1|}\vdash[cs'd]_l\dashv [b]_1$ is a normal $S$-diword,
and
$$
[u_1]_{|u_1|}\vdash[c\widetilde{s}'d]_l\dashv [b]_1\leq [u_1]_{|u_1|}\vdash\overline{([u_2a]_1\dashv s)}\dashv [b]_1
=\overline{([u]_n\dashv [asb]_m)}.
$$
If $a\widetilde{s}b< w$, $w\in X^+$, then $ua\widetilde{s}b< uw$ and $\overline{([u]_n\dashv [asb]_m)}\leq [ua\widetilde{s}b]_n<[uw]_n$.
 \ \ $\square$

 \ \

The lemma below follows from Lemma \ref{lnewa}, Definitions \ref{dcom} and \ref{dgsb} immediately.
\begin{lemma}\label{rr}
Let $S$ be a Gr\"{o}bner-Shirshov basis in $Di\langle X\rangle$, $f,g$ strong polynomials in $S$, $x\in X$, $[u]_n\in [X^+]_\omega$. Then the following statements hold.
\begin{enumerate}
\item[(i)]\
If  $w=\widetilde{f}b=a\widetilde{g}$ for some $a,b\in X^*$ such that $|\widetilde{f}|+|\widetilde{g}|>|w|$,
then
$$
[xfb]_1-[xag]_1\equiv 0 \ mod(S), \ \ [fbx]_{|wx|}-[agx]_{|wx|}\equiv 0 \ mod(S).
$$
Moreover,
\begin{eqnarray*}
[u]_n\dashv ([xfb]_1-[xag]_1)&\equiv& 0 \ mod(S,\ [uxw]_n),\\
([fbx]_{|w|+1}-[agx]_{|w|+1})\vdash [u]_n &\equiv& 0 \ mod(S, \ [wxu]_{|wx|+n}).
\end{eqnarray*}
\item[(ii)]\ If  $w=\widetilde{f}=a\widetilde{g}b$ for some $a,b\in X^*$,
then
$$
[xf]_1-[xagb]_1\equiv 0 \ mod(S), \ \ [fx]_{|wx|}-[agbx]_{|wx|}\equiv 0 \ mod(S).
$$
Moreover,
\begin{eqnarray*}
[u]_n\dashv([xf]_1-[xagb]_1)&\equiv& 0 \ mod(S,\ [uxw]_n),\\
([fx]_{|w|+1}-[agbx]_{|w|+1})\vdash [u]_n &\equiv& 0 \ mod(S, \ [wxu]_{|wx|+n}).
\end{eqnarray*}
\end{enumerate}
\end{lemma}

\begin{lemma}\label{l4}
Let $S$ be closed under the compositions of left and right multiplication.
Then for any $S$-diword $(asb)$, $(asb)$ has an expression:
$$
(asb)=\sum\alpha_i[a_is_ib_i]_{m_i},
$$
where each $\alpha_i\in \mathbf{k}, \ s_i\in S, \ a_i,b_i\in X^*$, and $\overline{[a_i s_i b_i]_{m_i}}\leq \overline{ (asb)}$.
\end{lemma}
\noindent{\bf Proof.} We may assume  that
$$
(asb)=[x_{i_1}\cdots x_{i_k}\cdots x_{i_n}]_m|_{x_{i_k}\mapsto s}.
$$
If $k=m$ or $s$ is strong, then $(asb)$ is a normal $S$-diword and the result holds.
Suppose that $k\neq m$ and $s$ is not strong. Then
$$
(asb)=([a]_{|a|}\vdash s)\vdash [b]_{m-|a|-1} \ \  \mathrm{or} \ \ (asb)=[a]_m\dashv (s\dashv [b]_1).
$$
Clearly, $[a]_{|a|}\vdash s$ and $s\dashv [b]_1$ are normal $S$-diwords.
In both cases, the result follows from Lemma \ref{lnewa}.  \ \ $\square$

\begin{lemma}\label{l00}
Let $S$ be a monic subset of $Di\langle X\rangle$.
Then for any nonzero polynomial $f\in Di\langle X\rangle$,
$$
f=\sum\alpha_{i}[u_{i}]_{n_i}+\sum\beta_{j}[a_js_jb_j]_{m_j},
$$
where each $[u_i]_{n_i}\in Irr(S),\ \alpha_i, \beta_j\in
\mathbf{k},\ a_j, b_j\in X^*,\ s_j\in S$, $[u_i]_{n_i}\leq \overline{f}$ and
$\overline{[a_js_jb_j]_{m_j}}\leq\overline{f}$.
\end{lemma}
\noindent{\bf Proof.}
Let $f=lc(f) \overline{f}+r\!_{_f}$. If $\overline{f}\in Irr(S)$,
then take $[u]_n=\overline{f}$ and $f_1=f-lc(f) [u]_n$. If $\overline{f}\notin Irr(S)$,
then $\overline{f}=\overline{[asb]_m}$ for some normal $S$-diword $[asb]_m$ and take $f_1=f-lc(f)[asb]_m$.
In both cases, we have $\overline{f_1}<\overline{f}$ and the result
follows from induction on $\overline{f}$.\ \ $\square$

\begin{lemma}\label{lkey}
Let $S$ be a Gr\"{o}bner-Shirshov basis in $Di\langle X\rangle$ and $[a_1s_1b_1]_{m_1},[a_2s_2b_2]_{m_2}$ normal $S$-diwords.
If $[w]_m=\overline{[a_1s_1b_1]_{m_1}}=\overline{[a_2s_2b_2]_{m_2}}$, then
$$
[a_1s_1b_1]_{m_1}-[a_2s_2b_2]_{m_2}\equiv 0 \ mod(S,[w]_m).
$$
\end{lemma}
\noindent{\bf Proof.} Since $[w]_m=\overline{[a_1s_1b_1]_{m_1}}=\overline{[a_2s_2b_2]_{m_2}}$,
it follows that $w=a_1\widetilde{s_1}b_1=a_2\widetilde{s_2}b_2$ and $m=m_1=m_2$.
Here we need consider three cases:
\begin{displaymath}
\begin{cases}
\text{ Case 1 $\widetilde{s_1}$ and $\widetilde{s_2}$ are mutually disjoint.}\\
\text{ Case 2 One of $\widetilde{s_1}$ and $\widetilde{s_2}$ is a subword of the other.}\\
\text{ Case 3 $\widetilde{s_1}$ and $\widetilde{s_2}$  have a nonempty intersection.}\\
\end{cases}
\end{displaymath}

For Case 1, we may assume that $\widetilde{s_1}$ is at the left of $\widetilde{s_2}$, i.e.
$b_1=a\widetilde{s_2}b_2$ and $a_2=a_1\widetilde{s_1}a$, here $a$ may be empty.
Then
$$
[a_2s_2b_2]_{m_1}-[a_1s_1b_1]_{m_1}=[a_1\widetilde{s_1}as_2b_2]_{m_1}-[a_1s_1a\widetilde{s_2}b_2]_{m_1}=: A.
$$
Let $s_1=\overline{s_1}+\sum \beta_i[u_i]_{n_i}$, $s_2=\overline{s_2}+\sum \beta'_j[v_j]_{l_j}$.
Here we have to discuss five cases:
\begin{displaymath}
\begin{cases}
\text{ Case 1.1 } \ 1\leq m_1\leq |a_1|, \ \text{which implies that $s_1,s_2$ are strong}.\\
\text{ Case 1.2 } \ m_1= |a_1|+p(\overline{s_1}), \ \text{which implies that $s_2$ is strong}.\\
\text{ Case 1.3 } \ |a_1\widetilde{s_1}|<m_1\leq |a_1\widetilde{s_1}a|, \ \text{which implies that $s_1,s_2$ are strong}.\\
\text{ Case 1.4 } \ m_1=|a_1\widetilde{s_1}a|+p(\overline{s_2}), \ \text{which implies that $s_1$ is strong}.\\
\text{ Case 1.5 } \ |a_1\widetilde{s_1}a\widetilde{s_2}|<m_1\leq |a_1\widetilde{s_1}a\widetilde{s_2}b_2|, \ \text{which implies that $s_1,s_2$ are strong}.
\end{cases}
\end{displaymath}
We first give the proof for Case 1.1. The same proof remains valid for Cases 1.3 and 1.5.
Since $1\leq m_1\leq |a_1|$, we have
\begin{eqnarray*}
A&=&[a_1]_{m_1}\dashv \overline{s_1}\dashv [a]_1\dashv s_2\dashv[b_2]_1
-[a_1]_{m_1}\dashv s_1\dashv [a]_1\dashv \overline{s_2}\dashv[b_2]_1\\
&=&-[a_1]_{m_1}\dashv (s_1-\overline{s_1})\dashv [a]_1\dashv s_2\dashv[b_2]_1
+[a_1]_{m_1}\dashv s_1\dashv [a]_1\dashv (s_2-\overline{s_2})\dashv[b_2]_1\\
&=&-\sum \beta_i[a_1]_{m_1}\dashv [u_i]_{n_i}\dashv [a]_1\dashv s_2\dashv[b_2]_1
+\sum \beta'_j[a_1]_{m_1}\dashv s_1\dashv [a]_1\dashv [v_j]_{l_j}\dashv[b_2]_1.
\end{eqnarray*}
As $s_1,s_2$ are strong we have
$$
[a_1]_{m_1}\dashv [u_i]_{n_i}\dashv [a]_1\dashv s_2\dashv[b_2]_1=[a_1u_ias_2b_2]_{m_1},
$$
$$
[a_1]_{m_1}\dashv s_1\dashv [a]_1\dashv [v_j]_{l_j}\dashv[b_2]_1=[a_1s_1av_jb_2]_{m_1},
$$
and
$u_i< \widetilde{s_1}, v_j< \widetilde{s_2}$. It follows that
$$
[a_1u_ia\widetilde{s_2}b_2]_{m_1}<[a_1\widetilde{s_1}a\widetilde{s_2}b_2]_{m_1}=[w]_{m_1}, \
[a_1\widetilde{s_1}av_jb_2]_{m_1}<[a_1\widetilde{s_1}a\widetilde{s_2}b_2]_{m_1}=[w]_{m_1}.
$$
We proceed to show Case 1.2. Similar proof applies to Case 1.4.
Since  $m_1= |a_1|+p(\overline{s_1})$, we have
\begin{eqnarray*}
A&=&[a_1]_{|a_1|}\vdash \overline{s_1}\dashv [a]_1\dashv s_2\dashv[b_2]_1
-[a_1]_{|a_1|}\vdash s_1\dashv [a]_1\dashv \overline{s_2}\dashv[b_2]_1\\
&=&-[a_1]_{|a_1|}\vdash (s_1-\overline{s_1})\dashv [a]_1\dashv s_2\dashv[b_2]_1
+[a_1]_{|a_1|}\vdash s_1\dashv [a]_1\dashv (s_2-\overline{s_2})\dashv[b_2]_1\\
&=&-\sum \beta_i[a_1]_{|a_1|}\vdash [u_i]_{n_i}\dashv [a]_1\dashv s_2\dashv[b_2]_1
+\sum \beta'_j[a_1]_{|a_1|}\vdash s_1\dashv [a]_1\dashv [v_j]_{l_j}\dashv[b_2]_1.
\end{eqnarray*}
It is clear that
$$
[a_1]_{|a_1|}\vdash s_1\dashv [a]_1\dashv [v_j]_{l_j}\dashv[b_2]_1=[a_1s_1av_jb_2]_{m_1}, \ [v_j]_{l_j}< \overline{s_2}.
$$
As $s_1$ is strong we also have
$$
[a_1]_{|a_1|}\vdash [u_i]_{n_i}\dashv [a]_1\dashv s_2\dashv[b_2]_1=[a_1u_ias_2b_2]_{|a_1|+n_i}, \ u_i< \widetilde{s_1}.
$$
It follows that
$$
[a_1\widetilde{s_1}av_jb_2]_{m_1}<[a_1\widetilde{s_1}a\widetilde{s_2}b_2]_{m_1}=[w]_{m_1}, \
[a_1u_ia\widetilde{s_2}b_2]_{|a_1|+n_i}<[a_1\widetilde{s_1}a\widetilde{s_2}b_2]_{m_1}=[w]_{m_1}.
$$

We now turn to Case 2, and may assume that $\widetilde{s_2}$ is a subword of $\widetilde{s_1}$, say,
$w'=\widetilde{s_1}=a\widetilde{s_2}b$. Then $a_2=a_1a$, $b_2=bb_1$ and
$$
[a_2s_2b_2]_{m_1}-[a_1s_1b_1]_{m_1}=[a_1as_2bb_1]_{m_1}-[a_1s_1b_1]_{m_1}=: B.
$$
Here we also should discuss five cases:
\begin{displaymath}
\begin{cases}
\text{ Case 2.1 } \ 1\leq m_1<|a_1|, \ \text{which implies that $s_1,s_2$ are strong}.\\
\text{ Case 2.2 } \ m_1=|a_1|, \ \text{which implies that $s_1,s_2$ are strong}.\\
\text{ Case 2.3 } \ m_1=|a_1|+p(\overline{s_1}),\
p(\overline{s_1})\in P[as_2b].\\
\text{ Case 2.4 } \ m_1=|a_1\widetilde{s_1}|+1, \ \text{which implies that $s_1,s_2$ are strong}.\\
\text{ Case 2.5 } \ |a_1\widetilde{s_1}|+1<m_1\leq |a_1\widetilde{s_1}b_1|, \ \text{which implies that $s_1,s_2$ are strong}.\\
\end{cases}
\end{displaymath}
In Cases 2.1 and 2.2, let $a_1=a'_1x$ and in Cases 2.4 and 2.5, let $b_1=yb'_1$, where $a'_1,b'_1\in X^*$, $x,y\in X$.
Then
\begin{numcases}{ }
\text{ Case 2.1 } \ B=-[a'_1]_{m_1}\dashv ([xs_1]_1-[xas_2b]_1)\dashv[b_1]_1. \label{eq1} \\
\text{ Case 2.2 } \ B=-[a'_1]_{|a'_1|}\vdash ([xs_1]_1-[xas_2b]_1)\dashv[b_1]_1. \label{eq2}\\
\text{ Case 2.3 } \ B=-[a_1]_{|a_1|}\vdash (s_1-[as_2b]_{p(\overline{s_1})})\dashv[b_1]_1. \label{eq3}\\
\text{ Case 2.4 } \ B=-[a_1]_{|a_1|}\vdash ([s_1y]_{|w'y|}-[as_2by]_{|w'y|})\dashv[b'_1]_1. \label{eq4}\\
\text{ Case 2.5 } \ B=-[a_1]_{|a_1|}\vdash ([s_1y]_{|w'y|}-[as_2by]_{|w'y|})\vdash[b'_1]_{m_1-|a_1\widetilde{s_1}|-1}. \label{eq5}
\end{numcases}
By Lemmas \ref{l0} and \ref{rr}, (\ref{eq1}) and (\ref{eq5}) are linear combinations of normal $S$-diwords
with leading monomials less than $[a'_1xw'b_1]_{m_1}=[a_1w'yb'_1]_{m_1}=[w]_{m_1}$.
Applying Lemmas \ref{l0}, \ref{rr} and using the fact that $S$ is a Gr\"{o}bner-Shirshov basis,
the same conclusion can be drawn for (\ref{eq2}), (\ref{eq3}) and (\ref{eq4}).

For Case 3, we may assume that $\widetilde{s_1}$ is at the left of $\widetilde{s_2}$, i.e.
$a_2=a_1a$, $b_1=bb_2$, and $w'=\widetilde{s_1}b=a\widetilde{s_2}$.
Then
$$
[a_2s_2b_2]_{m_1}-[a_1s_1b_1]_{m_1}=[a_1as_2b_2]_{m_1}-[a_1s_1bb_2]_{m_1}=:C.
$$
Here we continue to discuss five cases:
\begin{displaymath}
\begin{cases}
\text{ Case 3.1 } \ 1\leq m_1<|a_1|, \ \text{which implies that $s_1,s_2$ are strong}.\\
\text{ Case 3.2 } \ m_1=|a_1|, \ \text{which implies that $s_1,s_2$ are strong}.\\
\text{ Case 3.3 } \ m_1=|a_1|+m, m\in P([as_2])\cap P([s_1b]). \\
\text{ Case 3.4 } \ m_1=|a_1a\widetilde{s_2}|+1, \ \text{which implies that $s_1,s_2$ are strong}.\\
\text{ Case 3.5 } \ |a_1a\widetilde{s_2}|+1<m_1\leq |a_1a\widetilde{s_2}b_2|, \ \text{which implies that $s_1,s_2$ are strong}.\\
\end{cases}
\end{displaymath}
In Cases 3.1 and 3.2, let $a_1=a'_1x$ and in Cases 3.4 and 3.5, let $b_1=yb'_1$, where $a'_1,b'_1\in X^*$, $x,y\in X$.
Then
\begin{numcases}{ }
\text{ Case 3.1 } \ C=[a'_1]_{m_1}\dashv ([xas_2]_1-[xs_1b]_1)\dashv[b_2]_1. \label{eq6} \\
\text{ Case 3.2 } \ C=[a'_1]_{|a'_1|}\vdash ([xas_2]_1-[xs_1b]_1)\dashv[b_2]_1. \label{eq7}\\
\text{ Case 3.3 } \ C=[a_1]_{|a_1|}\vdash ([as_2]_m-[s_1b]_m)\dashv[b_2]_1. \label{eq8}\\
\text{ Case 3.4 } \ C=[a_1]_{|a_1|}\vdash ([as_2y]_{|w'y|}-[s_1by]_{|w'y|})\dashv[b'_2]_1. \label{eq9}\\
\text{ Case 3.5 } \ C=[a_1]_{|a_1|}\vdash ([as_2y]_{|w'y|}-[s_1by]_{|w'y|})\vdash[b'_2]_{m_1-|a_1a\widetilde{s_2}|-1}. \label{eq10}
\end{numcases}
By Lemmas \ref{l0} and \ref{rr}, (\ref{eq6}) and (\ref{eq10}) are linear combinations of normal $S$-diwords
with leading monomials less than $[a'_1xw'b_2]_{m_1}=[a_1w'yb'_2]_{m_1}=[w]_{m_1}$.
Applying Lemmas \ref{l0}, \ref{rr} and using the fact that $S$ is a Gr\"{o}bner-Shirshov basis,
the same conclusion can be drawn for (\ref{eq7}), (\ref{eq8}) and (\ref{eq9}).
\ \ $\square$

\begin{theorem}\label{cd}
(Composition-Diamond lemma for dialgebras) \ Let $S$ be a monic subset of $Di\langle X\rangle$ ,
$>$ a monomial-center ordering on $[X^+]_\omega$ and
$Id(S)$ the ideal of $Di\langle X\rangle$ generated by $S$. Then the following statements are equivalent.
\begin{enumerate}
\item[(i)] \ $S$ is a Gr\"{o}bner-Shirshov basis in $Di\langle X\rangle$.
\item[(ii)] \ $f\in Id(S)\Rightarrow \overline{f}=\overline{[asb]_m}$ for some normal $S$-diword $[asb]_m$.
\item[(iii)] \
$Irr(S)=\{[u]_n\in [X^+]_\omega\mid [u]_n\neq \overline{[asb]_m} \ \mbox{ for any normal } S\mbox{-diword}\ [asb]_m \}$
is a $\mathbf{k}$-basis of the quotient dialgebra $Di\langle X\mid  S\rangle:=Di\langle X\rangle/Id(S)$.
\end{enumerate}
\end{theorem}
\noindent{\bf Proof.} $(i)\Rightarrow (ii)$. Let $0\neq f\in Id(S)$. Then by
Lemma \ref{l4} $f$ has an expression
\begin{equation}\label{e4}
f=\sum \alpha_i[a_is_ib_i]_{m_i},
\end{equation}
where each $\alpha_i\in \mathbf{k}, \ a_i,b_i\in X^*, \ s_i\in S$. Write
$[w_i]_{m_i}=\overline{[a_is_ib_i]_{m_i}}=[a_i\widetilde{s_i}b_i]_{m_i}, i=1,2,\cdots$.
We may assume without loss of generality that
$$
[w_1]_{m_1}=[w_2]_{m_2}=\cdots=[w_l]_{m_l}>[w_{l+1}]_{m_{l+1}}\geq [w_{l+2}]_{m_{l+2}}\geq\cdots .
$$
The proof follows by induction on $([w_1]_{m_1},l)$. If $l=1$, then
$\overline{f}=\overline{[a_1s_1b_1]_{m_1}}=[a_1\widetilde{s_1}b_1]_{m_1}$ and
the result holds. Suppose that $l\geq 2$.
Then
$$
[w_1]_{m_1}=\overline{[a_1s_1b_1]_{m_1}}=\overline{[a_2s_2b_2]_{m_2}}.
$$
By Lemma \ref{lkey}, we can rewrite the first two summands of (\ref{e4}) in the form
\begin{eqnarray*}
\alpha_1[a_1s_1b_1]_{m_1}+\alpha_2[a_2s_2b_2]_{m_1}&=&(\alpha_1+\alpha_2)[a_1s_1b_1]_{m_1}+\alpha_2([a_2s_2b_2]_{m_1}-[a_1s_1b_1]_{m_1}) \\
&=&(\alpha_1+\alpha_2)[a_1s_1b_1]_{m_1}+\sum \alpha_2\beta_j[c_js'_jd_j]_{n_j},
\end{eqnarray*}
where each $[c_js'_jd_j]_{n_j}$ is a normal $S$-diword and $\overline{[c_js'_jd_j]_{n_j}}<[w_1]_{m_1}$.
Thus the result follows from induction on $([w_1]_{m_1},l)$.

$(ii)\Rightarrow (iii)$. By Lemma \ref{l00}, the set $Irr(S)$ generates
$Di\langle X\mid S\rangle$ as a linear space. On the other hand, suppose that
$h=\sum\alpha_i[u_i]_{l_i}=0$ in $Di\langle X\mid S\rangle$, where
each $\alpha_i\in \mathbf{k}$, $[u_i]_{l_i}\in {Irr(S)}$. This means that
$h\in Id(S)$.
Then all $\alpha_i$ must be equal to zero. Otherwise, $\overline{h}=[u_j]_{l_j}$ for some $j$ which contradicts $(ii)$.

$(iii)\Rightarrow(i)$. Suppose that $h$ is a composition
of elements of $S$. Clearly, $h\in Id(S)$. By Lemma \ref{l00},
$$
h=\sum_{i}\alpha_{i}[u_{i}]_{n_i}+\sum_{j}\beta_{j}[a_js_jb_j]_{m_j},
$$
where each $[u_i]_{n_i}\in Irr(S),\ \alpha_i, \beta_j\in
\mathbf{k},\ a_j, b_j\in X^*,\ s_j\in S$, and $[u_i]_{n_i}\leq\overline{h},\
\overline{[a_js_jb_j]_{m_j}}\leq\overline{h}$.
Then $\sum_{i}\alpha_{i}[u_{i}]_{n_i}\in Id(S)$.
By $(iii)$, we have $\alpha_{i}=0$ and
$h\equiv0 \ mod(S)$.
 \ \  $\square$

\begin{remark}
\emph{In \cite{Di}, a Composition-Diamond lemma for dialgebras is established and claims that}
$(i)\Rightarrow (iii)$, \emph{but not conversely. The reason is that the definitions of
a Gr\"{o}bner-Shirshov basis in} $Di\langle X\rangle$ \emph{are different, see Remark \ref{rcomp}}.
\end{remark}

\noindent{\bf Shirshov algorithm} If a monic subset $S \subset Di
\langle X\rangle$ is not a Gr\"{o}bner-Shirshov basis then one can add to $S$ all
nontrivial compositions. Continuing this process
repeatedly, we finally obtain a Gr\"{o}bner-Shirshov basis $S^{comp}$ that contains $S$
and generates the same ideal, $Id(S^{comp})=Id(S)$.

\begin{definition}
\emph{A Gr\"{o}bner-Shirshov basis} $S$ \emph{in} $Di\langle X\rangle$ \emph{is} minimal \emph{if for any} $s\in S$, $\overline{s}\in Irr(S \backslash \{s\})$.
\emph{A Gr\"{o}bner-Shirshov basis} $S$ \emph{in} $Di\langle X\rangle$ \emph{is} reduced \emph{if for any} $s\in S$, $supp(s)\subseteq Irr(S \backslash \{s\})$, \emph{where}
$$
supp(s):=\{[u_1]_{m_1},\dots,[u_n]_{m_n}\}
$$
\emph{if} $s=\alpha_1[u_1]_{m_1}+\dots+\alpha_n[u_n]_{m_n},\ 0\neq\alpha_i\in \mathbf{k},\ [u_i]_{m_i}\in [X^+]_\omega$.

\emph{Suppose} $I$ \emph{is an ideal of} $Di\langle X\rangle$ \emph{and} $I=Id(S)$. \emph{If} $S$ \emph{is a (reduced) Gr\"{o}bner-Shirshov basis in} $Di\langle X\rangle$, \emph{then we also call} $S$ \emph{is a (reduced) Gr\"{o}bner-Shirshov basis} for the ideal $I$ \emph{or} for the quotient dialgebra $Di\langle X\rangle/I$.
\end{definition}

For associative algebras and polynomial algebras, it is known that every ideal has a unique reduced Gr\"{o}bner-Shirshov basis.
This result is still true for dialgebras.

\begin{lemma}\label{lirr}
Let $I$ be an ideal of $Di\langle X\rangle$ and $S$ a Gr\"{o}bner-Shirshov basis for $I$. For any $T\subseteq S$,
if $Irr(T)=Irr(S)$ then $T$ is also a Gr\"{o}bner-Shirshov basis for $I$.
\end{lemma}
\noindent{\bf Proof.} For any $f\in I$, since $Irr(T)=Irr(S)$ and $S$ a Gr\"{o}bner-Shirshov basis for $I=Id(S)$, we have, by Theorem \ref{cd}, $ \overline{f}=\overline{[asb]_m}=\overline{[cgd]_m}$ for some $s\in S,\ g\in T,\ a,b,c,d\in X^*$. Thus,
$f_1=f-lc(f)[cgd]_m\in I$ and $\overline{f_1}<\overline{f}$. By induction on $\overline{f}$, $f$ is a linear combination of normal $T$-diwords, i.e. $f\in Id(T)$. This shows that $I= Id(T)$. Now the result follows from Theorem \ref{cd}. \ \ $\square$

Let $S$ be a subset of $Di\langle X\rangle$ and $[w]_m\in [X^+]_\omega$. We set
$$
\overline{S}:=\{\overline{s}\in [X^+]_\omega \mid s\in S\}, \ \ S^{[w]_m}:=\{s\in S\mid \overline{s}=[w]_m\}, \ \
S^{<{[w]_m}}:=\{s\in S\mid \overline{s}<[w]_m\}.
$$

\begin{theorem}
Let $I$ be an ideal of $Di\langle X\rangle$ and  $>$ a monomial-center ordering on $[X^+]_\omega$.
Then there is a unique reduced Gr\"{o}bner-Shirshov basis for $I$.
\end{theorem}
\noindent{\bf Proof.}
It is clear that there is a Gr\"{o}bner-Shirshov basis $S$ for $I$, for example, we may take $S=\{lc(f)^{-1}f\mid 0\neq f\in I\}$.
For each $[w]_m\in \overline{S}$, we choose a polynomial $f^{[w]_m}$ in $S$ such that $\overline{f^{[w]_m}}=[w]_m$. Write
$$
S_0=\{f^{[w]_m}\in S\mid [w]_m\in \overline{S}\}.
$$
Noting that $I\supseteq S\supseteq S_0$ and $\overline{I}=\overline{S}=\overline{S_0}$, we have  $Irr(S_0)=Irr(S)=[X^+]_\omega \backslash \overline{S}$.
By Lemma \ref{lirr}, $S_0$ is a Gr\"{o}bner-Shirshov basis for $I$.

Moreover, we may assume that for any $s\in S_0$,
\begin{equation}\label{eq99}
supp(s-\overline{s})\subseteq Irr(S_0)
\end{equation}
i.e.
$supp(s-\overline{s})\subseteq [X^+]_\omega \backslash \overline{S_0}$.
Indeed, if $supp(s-\overline{s})\cap \overline{S_0} \neq \emptyset$ for some $s\in S_0$, then set
$[u]_n=\max\{supp(s-\overline{s})\cap \overline{S_0}\}$ and there is an $f\in S_0$ such that $\overline{f}=[u]_n$.
Note that $\overline{s}>[u]_n=\overline{f}$ and $\overline{s-f}=\overline{s}$. Replace $s$ by $s-f$ in $S_0$.
Then $supp(s-f-\overline{s-f})\cap \overline{S_0} =\emptyset$  or $\max\{supp(s-f-\overline{s-f})\cap \overline{S_0}\}<[u]_n$.
Since $>$ is a well ordering on $[X^+]_\omega$, this process will terminate.

Note that for any $[w]_m\in \overline{S_0}$, there exists a unique $f\in S_0$ such that $[w]_m=\overline{f}$.
Set $\min\{\overline{S_0}\}=\overline{s_0}$ with $s_0\in S_0$. Define
$S_{_{\overline{s_0}}}:=\{s_0\}$.
Suppose that $f\in S_0,\ \overline{s_0}<\overline{f}$ and $S_{_{\overline{g}}}$ has been defined for any $g\in S_0$ with $\overline{g}<\overline{f}$. Define
\begin{displaymath} S_{_{\overline{ f}}}:=
\begin{cases}
S_{<\overline{ f}}
 \ \ \ \ \  \ \ \ \ \ \ \ \text{ if $\overline{ f}\not\in Irr(S_{<\overline{ f}})$,}\\
S_{<\overline{ f}}\cup \{f\}  \ \ \ \  \text{ if $\overline{ f}\in Irr(S_{<\bar f})$},
\end{cases}
\ \ \ \ \mbox{ where }\ \ S_{_{<\overline{ f}}}:=\bigcup_{ \overline{g}<  \overline{f},\ g\in S_0}S_{_{\overline{g}}}.
\end{displaymath}
Let
$$
S_1:=\bigcup_{f\in S_0} S_{_{\overline{ f}}}.
$$
Then  for any $f\in S_0,\ f\in S_1\Leftrightarrow \overline{ f}\in Irr(S_{<\bar f})\Leftrightarrow f\in S_{_{\overline{ f}}}$.

We first claim that $Irr(S_1)= Irr(S_0)$. Noting that $S_1\subseteq S_0$, it suffices to show that $Irr(S_1)\subseteq Irr(S_0)$. Assume that there is $[w]_m\in [X^+]_\omega$ such that $[w]_m\in Irr(S_1)$ and $[w]_m\notin Irr(S_0)$.
Since $\overline{S_0}=\overline{I}$, it follows that $[w]_m=\overline{f}$ for some $f\in S_0\backslash S_1$. If  $\overline{f}\in Irr(S_{<\overline{ f}})$ then $f\in S_{\overline{ f}}\subseteq S_1$, a contradiction. If  $\overline{f}\not\in Irr(S_{<\overline{ f}})$ then
$ \overline{f}=\overline{[asb]_m}$ for some $s\in S_{<\overline{ f}}\subseteq S_1, a,b\in X^*$.
This implies that  $\overline{ f}\not\in Irr(S_1)$, a contradiction. Therefore, $Irr(S_1)= Irr(S_0)$. Now by Lemma \ref{lirr}, $S_1$ is a  Gr\"{o}bner-Shirshov basis for $I$.

 If $f,g\in S_1, \ f\neq g,\ \overline{f}=\overline{[agb]_m}$, then $\overline{g}<\overline{f},\ g\in  S_{\overline{g}}\subseteq S_{<\overline{f}}$  which implies $ \overline{f}\not\in Irr(S_{<\overline{f}})$ and $ f\not\in S_1$, a contradiction. This shows that $S_1$ is a minimal Gr\"{o}bner-Shirshov basis for $I$.
By (\ref{eq99}), for any $s\in S_1$, $supp(s)\subseteq Irr(S_1 \backslash \{s\})$ which means $S_1$ is a reduced Gr\"{o}bner-Shirshov basis for $I$.

This shows that $I$ has a reduced Gr\"{o}bner-Shirshov basis $S_1$.

Suppose that $T$ is an arbitrary reduced Gr\"{o}bner-Shirshov basis for $I$. Let $\overline{s_0}=\min \overline{S_1}$ and
$\overline{r_0}=\min \overline{T}$, where $s_0\in S_1, r_0\in T$. By Theorem \ref{cd},
$\overline{s_0}=\overline{[ a'r'b']_p}\geq \overline{r'}\geq \overline{r_0}$ for some $r'\in T, a',b'\in X^*$.
Similarly, $\overline{r_0}\geq \overline{s_0}$. Then $\overline{r_0}=\overline{s_0}$. We say that $r_0=s_0$.
Otherwise, $0\neq r_0-s_0\in I$. We apply the above argument again, with replace $\overline{s_0}$ by $\overline{r_0-s_0}$,
to obtain that $\overline{r_0}> \overline{r_0-s_0}\geq\overline{r''}\geq \overline{r_0}$ for some $r''\in T$, a contradiction.
As both $T$ and $S_1$ are reduced Gr\"{o}bner-Shirshov bases, we have $S_1^{\overline{s_0}}=\{s_0\}=\{r_0\}=T^{\overline{r_0}}$.
Given any $[ w]_m\in \overline{S_1}\cup \overline{T}$ with $[ w]_m> \overline{r_0}$.
Assume that $S_1^{<[ w]_m}=T^{<[ w]_m}$.
To prove $T= S_1$, it is sufficient to show that $S_1^{[ w]_m} \subseteq T^{[ w]_m}$.
For any $s\in S_1^{[ w]_m}$, we can see that
$\overline{s}=\overline{[c'rd']_q}\geq \overline{r}$ for some $r\in T, c',d'\in X^*$.
Now, we claim that $[ w]_m=\overline{s}=\overline{r}$. Otherwise, $[ w]_m=\overline{s}>\overline{r}$.
Then $r\in T^{<[ w]_m}=S_1^{<[ w]_m}$ and $r\in S_1\backslash \{s\}$.
But $\overline{s}=\overline{[c'rd']_q}$, which contradicts the fact that $S_1$ is a reduced Gr\"{o}bner-Shirshov basis.
We next claim that $s=r\in T^{[ w]_m}$. If $s\neq r$, then $0\neq s-r\in I$.
By Theorem \ref{cd}, $\overline{s-r}=\overline{[ ar_1b]_n}=\overline{[ cs_1d]_n}$ for some
$r_1\in T, s_1\in S_1,a,b,c,d\in X^*$ with
$\overline{r_1}, \overline{s_1}\leq \overline{s-r}<\overline{s}=\overline{r}$.
This means that $s_1\in S_1\backslash \{s\}$ and $r_1\in T\backslash \{r\}$.
Noting that $\overline{s-r}\in supp(s)\cup supp(r)$, we may assume that $\overline{s-r}\in supp(s)$.
As $S_1$ is a reduced Gr\"{o}bner-Shirshov basis, we have $\overline{s-r}\in Irr(S_1\backslash \{s\})$,
which contradicts the fact that $\overline{s-r}=\overline{[ cs_1d]_n}$, where $s_1\in S_1\backslash \{s\}$.
Thus $s=r$. This shows that $S_1^{[ w]_m} \subseteq T^{[ w]_m}$.
\ \ $\square$

\begin{remark}
\emph{For associative algebras and polynomial algebras, it is known that every Gr\"{o}bner-Shirshov basis for an ideal can be reduced to a
reduced Gr\"{o}bner-Shirshov basis for the ideal.
Unfortunately, for dialgebras,  this is not the case.}
\end{remark}

The following example shows that generally, a Gr\"{o}bner-Shirshov basis $S$ in $Di\langle X\rangle$ may not be reduced to
a minimal Gr\"{o}bner-Shirshov basis for $I=Id(S)$.
\begin{example}
Let $X=\{x\}$, $Char\mathbf{k}\neq2,3$ and $S=\{f,g,h,p\}$,
where
\begin{eqnarray*}
&&f=[x^4]_4, \ \ g=[x^3]_3-\frac{1}{2}[x^3]_2-\frac{1}{2}[x^3]_1,\ \
h=[x^4]_3+[x^4]_2, \ \ p=[x^4]_2+\frac{1}{3}[x^4]_1.
\end{eqnarray*}
Then $S$
is a Gr\"{o}bner-Shirshov basis in $Di\langle X\rangle$
and $S$ can not be reduced to
a minimal Gr\"{o}bner-Shirshov basis for $I=Id(S)$.
\end{example}

\noindent{\bf Proof.} We first show that all compositions in $S$ are trivial.

1) Compositions of left (right) multiplication.

All possible compositions of left (right) multiplication are ones related to $g,h,p$.
By noting that for any $x^n\in X^+$, we have
$$
x\dashv g=0,  \ x\dashv h=2[x^5]_1=2x\dashv f \equiv 0 \ \mathrm{mod} (S), \ x\dashv p=\frac{4}{3}[x^5]_1=\frac{4}{3}x\dashv f \equiv 0 \ \mathrm{mod} (S);
$$
$$
g\vdash [x^n]_{n}=0, \ h\vdash [x^n]_{n}=2f\vdash [x^n]_{n}\equiv 0 \ \mathrm{mod} (S), \ p\vdash [x^n]_{n}=\frac{4}{3}f\vdash [x^n]_{n}\equiv 0 \ \mathrm{mod} (S).
$$

2) Compositions of inclusion and left (right) multiplicative inclusion.

We denote by, for example, ``$f\wedge g,\ [w]_m$" the composition of the polynomials of $f$ and  $g$ with ambiguity $[w]_m$.

By noting that in $S$,
\begin{eqnarray*}
&&f \wedge g,\ w=x^4, \
P(f)\cap P([xg])=\{4\},\ P(f)\cap P([gx])=\emptyset;\\
&&f \wedge h,\ w=x^4, \
P(f)\cap P(h)=\emptyset;\\
&&f \wedge p,\ w=x^4, \
P(f)\cap P(p)=\emptyset;\\
&&h \wedge g,\ w=x^4, \
P(h)\cap P([gx])=\{3\},\ P(h)\cap P([xg])=\emptyset;\\
&&h \wedge p,\ w=x^4, \
P(h)\cap P(p)=\emptyset;\\
&&p \wedge g,\ w=x^4, \
P(p)\cap P([xg])=\emptyset,\ P(p)\cap P([gx])=\emptyset,
\end{eqnarray*}
all possible of compositions of inclusion in $S$ are:
\begin{eqnarray*}
f\wedge g, [x^4]_4; \ \ \  h\wedge g, [x^4]_3.
\end{eqnarray*}
As $g,h,p$ are not strong, there is no composition of left (right) multiplicative inclusion.

For $f\wedge g$, $[w]_m=[x^4]_4$, we have
\begin{eqnarray*}
(f, g)_{[w]_m}=[x^4]_4-x\vdash ([x^3]_3-\frac{1}{2}[x^3]_2-\frac{1}{2}[x^3]_1)=\frac{1}{2}([x^4]_3+[x^4]_2)=\frac{1}{2}h\equiv 0 \ \mathrm{mod} (S).
\end{eqnarray*}

For $h\wedge g$, $[w]_m=[x^4]_3$, we have
\begin{eqnarray*}
(h, g)_{[w]_m}&=&[x^4]_3+[x^4]_2-([x^3]_3-\frac{1}{2}[x^3]_2-\frac{1}{2}[x^3]_1)\dashv x \\
&=&\frac{3}{2}[x^4]_2+\frac{1}{2}[x^4]_1=\frac{2}{3}p\equiv 0 \ \mathrm{mod} (S).
\end{eqnarray*}

3) Compositions of intersection and left (right) multiplicative intersection.

By noting that in $S$,
\begin{align*}
f \wedge f,\ w&=x^7, \
P([fx^3])\cap P([x^3f])=\{7\}; &
f \wedge f,\ w&=x^6, \
P([fx^2])\cap P([x^2f])=\{6\};\\
f \wedge f,\ w&=x^5, \
P([fx])\cap P([xf])=\{5\};&
f \wedge g,\ w&=x^6, \
P([fx^2])\cap P([x^3g])=\{6\};\\
f \wedge g,\ w&=x^5, \
P([fx])\cap P([x^2g])=\{5\}; &
f \wedge h,\ w&=x^7, \
P([fx^3])\cap P([x^3h])=\{6\};\\
f \wedge h,\ w&=x^6, \
P([fx^2])\cap P([x^2h])=\{5\};&
f \wedge h,\ w&=x^5, \
P([fx])\cap P([xh])=\{4\};\\
f \wedge p,\ w&=x^7, \
P([fx^3])\cap P([x^3p])=\{5\};&
f \wedge p,\ w&=x^6, \
P([fx^2])\cap P([x^2p])=\{4\};\\
f \wedge p,\ w&=x^5, \
P([fx])\cap P([xp])=\emptyset;&
g \wedge f,\ w&=x^6, \
P([gx^3])\cap P([x^2f])=\emptyset;\\
g \wedge f,\ w&=x^5, \
P([gx^2])\cap P([xf])=\emptyset;&
g \wedge g,\ w&=x^5, \
P([gx^2])\cap P([x^2g])=\emptyset;\\
g \wedge g,\ w&=x^4, \
P([gx])\cap P([xg])=\emptyset;&
g \wedge h,\ w&=x^6, \
P([gx^3])\cap P([x^2h])=\emptyset;\\
g \wedge h,\ w&=x^5, \
P([gx^2])\cap P([xh])=\emptyset;&
g \wedge p,\ w&=x^6, \
P([gx^3])\cap P([x^2p])=\emptyset;\\
g \wedge p,\ w&=x^5, \
P([gx^2])\cap P([xp])=\{3\};&
h \wedge f,\ w&=x^7, \
P([hx^3])\cap P([x^3f])=\{3\}; \\
h \wedge f,\ w&=x^6, \
P([hx^2])\cap P([x^2f])=\emptyset;&
h \wedge f,\ w&=x^5, \
P([hx])\cap P([xf])=\emptyset;\\
h \wedge g,\ w&=x^6, \
P([hx^2])\cap P([x^3g])=\emptyset;&
h \wedge g,\ w&=x^5, \
P([hx])\cap P([x^2g])=\emptyset;\\
h \wedge h,\ w&=x^7, \
P([hx^3])\cap P([x^3h])=\emptyset;&
h \wedge h,\ w&=x^6, \
P([hx^2])\cap P([x^2h])=\emptyset;\\
h \wedge h,\ w&=x^5, \
P([hx])\cap P([xh])=\emptyset;&
h \wedge p,\ w&=x^7, \
P([hx^3])\cap P([x^3p])=\emptyset;\\
h \wedge p,\ w&=x^6, \
P([hx^2])\cap P([x^2p])=\emptyset;&
h \wedge p,\ w&=x^5, \
P([hx])\cap P([xp])=\{3\};\\
p \wedge f,\ w&=x^7, \
P([px^3])\cap P([x^3f])=\{2\}; &
p \wedge f,\ w&=x^6, \
P([px^2])\cap P([x^2f])=\{2\};\\
f \wedge f,\ w&=x^5, \
P([px])\cap P([xf])=\emptyset;&
p \wedge g,\ w&=x^6, \
P([px^2])\cap P([x^3g])=\emptyset;\\
p \wedge g,\ w&=x^5, \
P([px])\cap P([x^2g])=\emptyset; &
p \wedge h,\ w&=x^7, \
P([px^3])\cap P([x^3h])=\emptyset;\\
p \wedge h,\ w&=x^6, \
P([px^2])\cap P([x^2h])=\emptyset;&
p \wedge h,\ w&=x^5, \
P([px])\cap P([xh])=\emptyset;\\
p \wedge p,\ w&=x^7, \
P([px^3])\cap P([x^3p])=\emptyset;&
p \wedge p,\ w&=x^6, \
P([px^2])\cap P([x^2p])=\emptyset;\\
p \wedge p,\ w&=x^5, \
P([px])\cap P([xp])=\emptyset,&
\end{align*}
all possible ambiguities $[w]_m$ of compositions of intersection are:
\begin{align*}
f&\wedge f, [x^7]_7,[x^6]_6,[x^5]_5; & f&\wedge g, [x^6]_6,[x^5]_5; &  f&\wedge h, [x^7]_6,[x^6]_5,[x^5]_4; & f&\wedge p, [x^7]_5,[x^6]_4;\\
g&\wedge p, [x^5]_3; & h&\wedge f, [x^7]_3; & h&\wedge p, [x^5]_3; & p&\wedge f, [x^7]_2,[x^6]_2.
\end{align*}
As $g,h,p$ are not strong, there is no composition of left (right) multiplicative intersection.
It is easy to see that all the compositions of intersection are trivial modulo $S$.

This shows that $S$ is a Gr\"{o}bner-Shirshov basis.

Note that in $S$, $\overline{f}=\overline{[xg]_4}$, $[xg]_4$ is a normal $g$-diword, and $f-[xg]_4=\frac{1}{2}h$. Since $2[x^5]_1=x\dashv h$ is nontrivial modulo $\{g,h,p\}$, $\{g,h,p\}$ is not a Gr\"{o}bner-Shirshov basis. This implies that we can not drop $f$ from $S$, i.e.  $S$ can not be reduced to
a minimal Gr\"{o}bner-Shirshov basis.
\ \ $\square$

\section{Gr\"{o}bner-Shirshov bases for  dirings}

In this section, by similar proofs of the above section, we introduce Gr\"{o}bner-Shirshov bases for dirings, which may find an $R$-basis for some disemigroup-dirings over an associative ring $R$.

\begin{definition}(\cite{Lo99})\label{l2}
\emph{A} disemigroup \emph{is a set} $D$ \emph{equipped with two maps}
$$
\vdash \ : D\times D\rightarrow D,
\ \
\dashv \ : D\times D\rightarrow D,
$$
\emph{where} $\vdash$ \emph{and} $\dashv$ \emph{are associative and satisfy the identities (\ref{eq00})}.

\end{definition}

Note that in \cite{Lo99,Zhuchok,Zhuchok17,Zhu15}, such a disemigroup in the above definition is called a dimonoid.

It is well known from \cite{Lo99} that $([X^+]_\omega,\vdash,\dashv)$ is the free disemigroup generated by $X$, where
for any $[u]_m, [v]_n \in [X^+]_\omega$,
$$
[u]_m\vdash [v]_n=[uv]_{|u|+n}, \ \ \ [u]_m\dashv [v]_n=[uv]_m.
$$

Let us denote
$$
Disgp\langle X\rangle:=([X^+]_\omega,\vdash,\dashv)
$$
the free disemigroup generated by $X$.

Throughout this section, $R$ is an associative ring with unit.
\begin{definition}
\emph{A} diring \emph{is a quaternary} $(T,+, \vdash, \dashv)$  \emph{such that both}
 $(T,+, \vdash)$ \emph{and} $(T,+,\dashv)$ \emph{are associative rings with the identities (\ref{eq00})
in}  $T$.
\end{definition}

\begin{definition}
\emph{Let} $(D, \vdash, \dashv)$ \emph{be a disemigroup}, $R$ \emph{an associative ring with unit and} $T$ \emph{the free left} $R$-\emph{module with} $R$-\emph{basis} $D$. \emph{Then}
$(T,+, \vdash, \dashv)$ \emph{is a diring with a natural way: for any} $f=\sum_{i} r_iu_i,\ g=\sum_{j} r_j'v_j\in T,\ r_i,r_j'\in R,\ u_i, v_j\in D$,
$$
f\vdash g:=\sum_{i,j} r_ir_j'u_i\vdash v_j,\ \ \ \ f\dashv g:=\sum_{i,j} r_ir_j'u_i\dashv v_j.
$$
\emph{Such a diring, denoted by $Di_R(D)$, is called a} disemigroup-diring \emph{of} $D$ \emph{over} $R$.

\emph{We denote by} $Di_R\langle X\rangle$ \emph{the disemigroup-diring of} $Disgp\langle X\rangle$ \emph{over} $R$
\emph{which is also called the}  free diring \emph{over} $R$ \emph{generated by} $X$.
\emph{In particular}, $Di_{\mathbf{k}}\langle X\rangle=Di\langle X\rangle$ \emph{is the free dialgebra} \emph{generated by} $X$ \emph{when} $\mathbf{k}$ \emph{is a field}.

\emph{An}  ideal $I$ \emph{of} $Di_R\langle X\rangle$ \emph{is an} $R$-\emph{submodule of} $Di_R\langle X\rangle$ \emph{such that}
$f\vdash g, f\dashv g,\ g\vdash f, g\dashv f \in I$ \emph{for any} $f\in Di_R\langle X\rangle$ \emph{and} $g\in I$.
\end{definition}

As  same as the proof of Theorem \ref{cd}, we have the following Composition-Diamond lemma for dirings.

\begin{theorem}\label{cd-}
(Composition-Diamond lemma for dirings) \ Let $S$ be a monic subset of $Di_R\langle X\rangle$,
$>$ a monomial-center ordering on $[X^+]_\omega$ and
$Id(S)$ the ideal of $Di_R\langle X\rangle$ generated by $S$. Then the following statements are equivalent.
\begin{enumerate}
\item[(i)] \ $S$ is a Gr\"{o}bner-Shirshov basis in $Di_R\langle X\rangle$.
\item[(ii)] \ $f\in Id(S)\Rightarrow \overline{f}=\overline{[asb]_m}$ for some normal $S$-diword $[asb]_m$.
\item[(iii)] \
$Irr(S)=\{[u]_n\in [X^+]_\omega\mid [u]_n\neq \overline{[asb]_m}$ for any normal $S$-diword $[asb]_m$ \}
is an $R$-basis of the quotient diring $Di_R\langle X\mid S\rangle:=Di_R\langle X\rangle/Id(S)$, i.e. $Di_R\langle X\mid S\rangle$ is a free $R$-module with $R$-basis $Irr(S)$.
\end{enumerate}
\end{theorem}

\begin{remark}
\emph{Shirshov algorithm does not work generally in} $Di_R\langle X\rangle$.
\end{remark}

\section{Applications}

In this section, by using our Theorem  \ref{cd}, we give a method to find normal forms of elements of an arbitrary disemigroup, in particular, we give normal forms of elements of free commutative disemigroups, free abelian disemigroups and free left (right) commutative disemigroups.

\subsection{Normal forms of disemigroups}\label{section3.2.1}

For an arbitrary disemigroup $D$, $D$ has an expression
$$
D=Disgp\langle X\mid S\rangle:=[X^+]_\omega/\rho(S)
$$
 for some set $X$ and $S\subseteq [X^+]_\omega\times[X^+]_\omega$, where $\rho(S)$ is the congruence on $([X^+]_\omega,\vdash,\dashv)$ generated by $S$.

It is natural to ask how to find normal forms of elements of  disemigroup  $Disgp\langle X\mid S\rangle$?

Let $>$ be a monomial-center ordering on $[X^+]_\omega$ and $S=\{([u_i]_{m_i},[v_i]_{n_i})\mid [u_i]_{m_i}>[v_i]_{n_i},\ i\in I\}$. Consider the dialgebra $Di\langle X\mid S\rangle$, where $S=\{[u_i]_{m_i}-[v_i]_{n_i}\mid  i\in I\}$.
By Shirshov algorithm, we have a Gr\"{o}bner-Shirshov basis $S^{comp}$ in $Di\langle X\rangle$ and $Id(S^{comp})=Id(S)$. It is clear that each element in $S^{comp}$ is of the form $[u]_m-[v]_n$, where $[u]_m>[v]_n,\ [u]_m,[v]_n\in [X^+]_\omega$. Let
\begin{eqnarray*}
\sigma&:&  \ Di\langle X\mid S\rangle\rightarrow Di_{\mathbf{k}}([X^+]_\omega/\rho(S)), \\
&&\sum\alpha_i[u_i]_{m_i}+Id(S)\mapsto \sum\alpha_i [u_i]_{m_i}\rho(S),\ \ \ \ \alpha_i\in \mathbf{k},\  [u_i]_{m_i}\in [X^+]_\omega.
\end{eqnarray*}
Then  $\sigma$ is obviously a  dialgebra isomorphism. Since $Irr(S^{comp})$ is a  linear basis of $
Di\langle X\mid S\rangle$, we have $\sigma(Irr(S^{comp}))$ is a linear basis of
$Di_{\mathbf{k}}([X^+]_\omega/\rho(S))$ which shows that $Irr(S^{comp})$ is exactly  normal
forms of elements of the disemigroup $Disgp\langle X\mid S\rangle$.

Therefore, we have the following theorem.

\begin{theorem}\label{tndsg}
Let $>$ be a monomial-center ordering on $[X^+]_\omega$ and $D=Disgp\langle X\mid S\rangle$, where $S=\{([u_i]_{m_i},[v_i]_{n_i})\mid [u_i]_{m_i}>[v_i]_{n_i},\ i\in I\}$ is a subset of $[X^+]_\omega\times[X^+]_\omega$. Then  $Irr(S^{comp})$ is a set of normal forms of elements of  the disemigroup
$Disgp\langle X\mid S\rangle$.
\end{theorem}

From now on, let $>$ be the deg-lex-center ordering on $[X^+]_\omega$, where $X$ is a well-ordered set.

\subsection{Normal forms of free commutative disemigroups}

The commutative disemigroups are introduced and the free commutative disemigroup generated by a set is constructed by \cite{Zhuchok}. In this subsection, we give another approach to normal forms of elements of a free commutative  disemigroup.

\begin{definition} (\cite{Zhuchok})
\emph{A disemigroup} $(D,\vdash,\dashv)$ \emph{is} commutative \emph{if both} $\vdash$ \emph{and} $\dashv$ \emph{are commutative}.
\end{definition}

Let $Di[X]$ be the free commutative dialgebra generated by a set $X$ and $T$ be the subset of $Di\langle X\rangle$
consisting of the following polynomials:
$$
[u]_m\vdash [v]_n-[v]_n\vdash [u]_m, \ \ \ \ \ \ [u]_m\dashv [v]_n-[v]_n\dashv [u]_m,
$$
where $[u]_m,[v]_n\in [X^+]_\omega$.
Then $Di[X]=Di\langle X\mid T\rangle$ and $Disgp[X]=Disgp\langle X\mid T\rangle$ is the free commutative  disemigroup generated by $X$.

Let $X=\{x_i\mid i\in I\}$ be a total-ordered set,
$$
\lfloor X^+\rfloor:= \{\lfloor x_{i_1}x_{i_2}\cdots  x_{i_n}\rfloor\mid  i_1,\dots,i_n\in I, i_1\leq i_2\leq \cdots \leq i_n, n\in \mathbb{Z}^+\}
$$
the set of all nonempty commutative associative words on $X$ and
$$
\lfloor X^+\rfloor_\omega:=\{\lfloor u\rfloor_m\mid \lfloor u\rfloor \in \lfloor X^+\rfloor, m\in \mathbb{Z}^+, 1\leq m \leq |u|\}
$$
the set of all commutative normal diwords on $X$.
For $ u\in X^+,\ [u]_m$ is called an \textit{associative diword}, while $\lfloor u\rfloor_m$ is called a \textit{commutative diword}. For example, if $u=x_2x_1x_2x_1\in X^+,\ x_1<x_2$, then $\lfloor u\rfloor=\lfloor x_1x_1 x_2x_2\rfloor,\ [u]_3=x_2x_1 \dot{x_2}x_1,\  \lfloor u\rfloor_3=\lfloor x_1x_1x_2x_2\rfloor_3=x_1x_1 \dot{x_2}x_2$.

\begin{proposition}\label{pfcd}
Let $X=\{x_i\mid i\in I\}$ be a well-ordered set. Then
\begin{enumerate}
\item[(i)] \ $Di[X]=Di\langle X\mid S\rangle$, where $S$ consists of the following
polynomials:
\begin{eqnarray*}
[u]_m-\lfloor u\rfloor_m,  \   ([u]_m\in [X^+]_\omega, |u|=2), \ \ \  [v]_n-\lfloor v\rfloor_1,  \  ([v]_n\in [X^+]_\omega, |v|\geq 3).
\end{eqnarray*}
\item[(ii)] \ $S$ is a Gr\"{o}bner-Shirshov basis in $Di\langle X\rangle$.
\item[(iii)] \ The set
$\lfloor X^+ \rfloor_{_1} \cup \lfloor X^+ \rfloor_{_{2-2}}$
is a $\mathbf{k}$-basis of the free commutative dialgebra $Di[X]$,
where
$$
\lfloor X^+ \rfloor_{_1}:=\{\lfloor v \rfloor_1 \mid \lfloor v\rfloor\in \lfloor X^+ \rfloor\}\ \mbox{ and } \  \lfloor X^+ \rfloor_{_{2-2}}:=\{\lfloor u \rfloor_2 \mid \lfloor u\rfloor\in \lfloor X^+ \rfloor, |u|=2\}.
$$
\end{enumerate}
\end{proposition}
\noindent{\bf Proof.}
$(i)$ We only need to prove that the polynomials in $S$ are trivial modulo $T$ and the polynomials in $T$ are trivial modulo $S$.
It is clear that
$$
[x_ix_j]_2-\lfloor x_jx_i\rfloor_2 \equiv 0 \ \ \mathrm{mod} (T), \ \ \ \ \  [v]_1-\lfloor v\rfloor_1 \equiv 0 \ \ \mathrm{mod} (T),
$$
where $x_i,x_j\in X, v\in X^+, |v|\geq 2$. Suppose that $v=x_{j_1}\cdots x_{j_n} \cdots x_{j_l}\in X^+$, $n\geq 2,l>2$.
If $n<l$, then $v=v_1x_{j_n}v_2$ for some $v_1,v_2\in X^+$ and
\begin{eqnarray*}
[v]_n-\lfloor v\rfloor_1&=&[v_1x_{j_n}]_n\dashv[v_2]_1-\lfloor v\rfloor_1
\equiv [v_2]_1\dashv[v_1x_{j_n}]_n-\lfloor v\rfloor_1 \\
&\equiv& [v_2v_1x_{j_n}]_1-\lfloor v\rfloor_1
\equiv \lfloor v_2v_1x_{j_n}\rfloor_1-\lfloor v\rfloor_1\equiv 0 \ \ \mathrm{mod} (T).
\end{eqnarray*}
If $n=l$, then $v=x_{j_1}v'x_{j_n}$ for some $v'\in X^+$ and
\begin{eqnarray*}
[v]_n-\lfloor v\rfloor_1&=&[x_{j_1}v']_1\vdash x_{j_n}-\lfloor v\rfloor_1
\equiv x_{j_n}\vdash[x_{j_1}v']_1-\lfloor v\rfloor_1 \\
&\equiv &[x_{j_n}x_{j_1}v']_2-\lfloor v\rfloor_1
\equiv \lfloor x_{j_n}x_{j_1}v'\rfloor_1-\lfloor v\rfloor_1 \equiv 0 \ \ \mathrm{mod} (T).
\end{eqnarray*}

It is easily seen that
$$
x\vdash y-y\vdash x \equiv 0 \ \ \mathrm{mod} (S),\ \ \ \ x\dashv y-y\dashv x \equiv 0 \ \ \mathrm{mod} (S),
$$
where $x,y\in X$.
Suppose that $[u]_m,[v]_n\in [X^+]_\omega$ with $|uv|>2$.
\begin{align*}
[u]_m\vdash [v]_n-[v]_n\vdash [u]_m&=[uv]_{|u|+n}-[vu]_{|v|+m} \equiv \lfloor uv\rfloor_1-\lfloor vu\rfloor_1 \equiv 0 \ \ \mathrm{mod} (S),\\
[u]_m\dashv [v]_n-[v]_n\dashv [u]_m&=[uv]_{m}-[vu]_{n} \equiv \lfloor uv\rfloor_1-\lfloor vu\rfloor_1 \equiv 0 \ \ \mathrm{mod} (S).
\end{align*}

$(ii)$ It is easy to check that all possible compositions of left (right) multiplication in $S$ are equal to zero.
For any composition of $(f,g)_{[w]_m}$ in $S$, note that $-r\!_{_f},-r\!_{_g}\in [X^+]_\omega$, $|w|\geq 3$,
$[w]_m=\overline{[afb]_m}=\overline{[cgd]_m}$ and
$\lfloor w\rfloor_1=\lfloor a\widetilde{r\!_{_f}}b\rfloor_1=\lfloor c\widetilde{r\!_{_g}}d\rfloor_1$,
where $f=\overline{f}+r\!_{_f},g=\overline{g}+r\!_{_g}$, $a,b,c,d\in X^*$. It follows that
\begin{eqnarray*}
(f,g)_{[w]_m}=[afb]_m-[cgd]_m=-[a\widetilde{r\!_{_f}}b]_{m_1}+[c\widetilde{r\!_{_g}}d]_{m_2}
\equiv -\lfloor a\widetilde{r\!_{_f}}b\rfloor_{1}+\lfloor c\widetilde{r\!_{_g}}d\rfloor_{1}\equiv 0 \ \  \mathrm{mod} (S).
\end{eqnarray*}
Then all the compositions in $S$ are trivial. We have proved $(ii)$.

$(iii)$ This part follows from  Theorem \ref{cd}.
\ \ $\square$

From Theorem \ref{cd},  Lemma \ref{lirr} and Proposition \ref{pfcd}, it follows that
\begin{corollary}
Let $W$ be a set consisting of the following
polynomials:
\begin{align*}
&[x_ix_j]_2-[x_jx_i]_2,\ \ [x_ix_j]_1-[x_jx_i]_1, & (&i,j\in I,\ i>j),\\
&[x_ix_jx_k]_2-[x_ix_jx_k]_1, \ \ [x_ix_jx_k]_3-[x_ix_jx_k]_1, & (&i,j,k\in I, \ i\leq j \leq k).
\end{align*}
Then
$W$ is the reduced Gr\"{o}bner-Shirshov basis for the free commutative dialgebra $Di[X]$.
\end{corollary}

From Theorem \ref{tndsg} and Proposition \ref{pfcd}, it follows that

\begin{corollary}(\cite[Theorem 3]{Zhuchok})
$Disgp[X]=( \lfloor X^+ \rfloor_{_1}\cup \lfloor X^+ \rfloor_{_{2-2}}, \ \vdash, \dashv \ )$ is the free commutative disemigroup generated by $X$,
where the operations $\vdash$ and $\dashv$ are as follows: for any
$x,x'\in X,\ \lfloor u\rfloor_{p_1}, \lfloor v\rfloor_{p_2}\in \lfloor X^+ \rfloor_{_1}\cup \lfloor X^+ \rfloor_{_{2-2}}$ with  $|u||v|>1$,
\begin{align*}
\lfloor v\rfloor_{p_2}\vdash \lfloor u\rfloor_{p_1}
&=\lfloor u\rfloor_{p_1}\vdash \lfloor v\rfloor_{p_2}=\lfloor u\rfloor_{p_1}\dashv \lfloor v\rfloor_{p_2}
=\lfloor v\rfloor_{p_2}\dashv \lfloor u\rfloor_{p_1}=\lfloor uv\rfloor_1,\\
x\dashv x'&=x'\dashv x=\lfloor xx'\rfloor_1, \\
x\vdash x'&=x'\vdash x=\lfloor xx'\rfloor_2.
\end{align*}
\end{corollary}

\subsection{Normal forms of free abelian disemigroups}

The concept of abelian disemigroups is introduced and the free abelian disemigroup generated by a set is constructed by \cite{Zhu15}. In this subsection, we give another approach to normal forms of elements of a free abelian disemigroup.

\begin{definition} (\cite{Zhu15})
\emph{A disemigroup} $(D,\vdash,\dashv)$ \emph{is} abelian \emph{if} $ a\vdash b= b\dashv a$ \emph{for all} $a,b\in D$.

\end{definition}

Let $X$ be an arbitrary set and $T$ the subset of $[X^+]_\omega\times[X^+]_\omega$
consisting of the following:
$$
([u]_m\vdash [v]_n,[v]_n\dashv [u]_m),
$$
where $[u]_m,[v]_n\in [X^+]_\omega$. Then $Disgp\langle X\mid T\rangle$ is the free abelian disemigroup generated by $X$.

Let $X=\{x_i\mid i\in I\}$ be a total-ordered set.
Suppose that $u=x_{j_1}x_{j_2}\cdots x_{j_n}\in X^+$ and $\lfloor u\rfloor=\lfloor x_{i_1}x_{i_2}\cdots x_{i_n}\rfloor$,
where $x_{i_1},x_{i_2},\cdots, x_{i_n}$ is the reordering of
$x_{j_1},x_{j_2},\cdots, x_{j_n}$ such that $ x_{i_1}\leq x_{i_2}\leq\cdots\leq x_{i_n}$.
Define
\begin{align*}
cont(u):&= \{x\in X\mid x=x_{j_t} \ \mbox{ for some }\ 1\leq t\leq n\}, \ \ \ L(u):=\{1,2,\cdots,n\}. \\
\rho_{_u}:& \ L(u)\rightarrow cont(u),\ \ \ \ \ \ \ \ \ m \mapsto x_{j_m}. \\
\lambda_{_{\lfloor u \rfloor}}:& \ cont(u)\rightarrow L(u), \ \ \ \ \ \ \ \ \ \ x\mapsto 1 \  \ \ \emph{if} \  \ \ \ x=x_{i_1}, \\
& \ \ \ \ \ \ \ \ \ \ \ \ \ \ \  \ \ \ \ \  \ \ \ \ \ \ \ \ \ \ \ \ \ x \mapsto l \  \ \ \emph{if}  \  \ \ \
 x=x_{i_l}, \ x_{i_l}>x_{i_{l-1}},\  l>1.\\
\tau_{_u}=\lambda_{_{\lfloor u \rfloor}}\rho_{_u}:& \ L(u)\rightarrow L(u).
\end{align*}
For example, if $u=x_1x_2x_1x_2$ with $x_1<x_2$, then $\lfloor u\rfloor=\lfloor x_1x_1x_2x_2\rfloor$, $cont(u)=\{x_1,x_2\}$, $L(u)=\{1,2,3,4\}$,
$\rho_{_u}(2)=x_2$, $\lambda_{_{\lfloor u \rfloor}}(x_2)=3$, $\tau_{_u}(2)=3$.

For any $u,v\in X^+$, it is easy to check that $\tau_{_{\lfloor u\rfloor} }(m)\leq m$ for all $m\in L(u)$
and $\tau_{_{uv} }(|u|+n)=\tau_{_{vu} }(n)$ for all $n\in L(v)$.

\begin{proposition}\label{pfad}
Let $X=\{x_i\mid i\in I\}$ be a well-ordered set,
$T$ the subset of $Di\langle X\rangle$
consisting of the following polynomials:
$[u]_m\vdash [v]_n-[v]_n\dashv [u]_m$, where $[u]_m,[v]_n\in [X^+]_\omega$. Then
\begin{enumerate}
\item[(i)] \ $Di\langle X\mid T\rangle=Di\langle X\mid S\rangle$, where $S$ consists of the following
polynomials:
\begin{align*}
[u]_m-\lfloor u\rfloor_{\tau_{_u}(m)} \ \ \ \ ([u]_m\in [X^+]_\omega,\ |u|\geq 2).
\end{align*}
\item[(ii)] \ $S$ is a Gr\"{o}bner-Shirshov basis in $Di\langle X\rangle$.
\item[(iii)] \ The set
\begin{align*}
\{\lfloor x_{i_1}\cdots x_{i_{t-1}}x_{i_t}x_{i_{t+1}}\cdots x_{i_l}\rfloor_t\mid &
x_{i_1}\leq\cdots\leq x_{i_{t-1}}<x_{i_t}\leq x_{i_{t+1}}\leq\cdots \leq x_{i_l}, \\
&x_{i_j}\in X,\ 1\leq j\leq l,\ \  l,t\in \mathbb{Z}^+, \ t\leq l
\}
\end{align*}
is a $\mathbf{k}$-basis of the free abelian  dialgebra $Di\langle X\mid T\rangle$.
\end{enumerate}
\end{proposition}
\noindent{\bf Proof.} $(i)$
We only need to prove that the polynomials in $S$ are trivial modulo $T$ and the polynomials in $T$ are trivial modulo $S$.
\begin{align*}
[u]_m\vdash [v]_n-[v]_n\dashv [u]_m=[uv]_{|u|+n}-[vu]_n\equiv \lfloor uv\rfloor_{\tau_{_{uv}}(|u|+n)}-\lfloor vu\rfloor_{\tau_{_{vu} }(n)}
 \equiv 0\ \ \mathrm{mod} (S).
\end{align*}
On the other hand, it is easy to see that $[u]_m-\lfloor u\rfloor_{\tau_{_u}(m)}\equiv 0 \ \ \mathrm{mod} (T)$.

$(ii)$ It is easy to check that all possible compositions of left (right) multiplication in $S$ are equal to zero.
For any composition of $(f,g)_{[w]_m}$ in $S$, note that $-r\!_{_f},-r\!_{_g}\in [X^+]_\omega$, $|w|\geq 3$,
$[w]_m=\overline{[afb]_m}=\overline{[cgd]_m}$ and
$\lfloor w\rfloor=\lfloor a\widetilde{r\!_{_f}}b\rfloor=\lfloor c\widetilde{r\!_{_g}}d\rfloor$,
where $f=\overline{f}+r\!_{_f},g=\overline{g}+r\!_{_g}$, $a,b,c,d\in X^*$. It follows that
\begin{eqnarray*}
(f,g)_{[w]_m}=[afb]_m-[cgd]_m=-[a\widetilde{r\!_{_f}}b]_{m_1}+[c\widetilde{r\!_{_g}}d]_{m_2}.
\end{eqnarray*}
From the definition of composition in $S$ we conclude that
$\rho_{_w}(m)=\rho_{_{a\widetilde{r\!_{_f}}b}}(m_1)=\rho_{_{c\widetilde{r\!_{_g}}d}}(m_2)$.
Thus $\tau_{_{a\widetilde{r\!_{_f}}b}}(m_1)=\lambda_{_{\lfloor a\widetilde{r\!_{_f}}b\rfloor}}\rho_{_{a\widetilde{r\!_{_f}}b}}(m_1)
=\lambda_{_{\lfloor c\widetilde{r\!_{_g}}d\rfloor}}\rho_{_{c\widetilde{r\!_{_g}}d}}(m_2)=\tau_{_{c\widetilde{r\!_{_g}}d}}(m_2)$
and
\begin{eqnarray*}
-[a\widetilde{r\!_{_f}}b]_{m_1}+[c\widetilde{r\!_{_g}}d]_{m_2}\equiv -\lfloor a\widetilde{r\!_{_f}}b\rfloor_{\tau_{_{a\widetilde{r\!_{_f}}b}}(m_1)}
+\lfloor c\widetilde{r\!_{_g}}d\rfloor_{\tau_{_{c\widetilde{r\!_{_g}}d}}(m_2)}
\equiv 0 \ \ \mathrm{mod} (S).
\end{eqnarray*}
Then all the compositions in $S$ are trivial. We have proved $(ii)$.

$(iii)$ This part follows from  Theorem \ref{cd}.
\ \ $\square$

From Theorem \ref{cd},  Lemma \ref{lirr} and Proposition \ref{pfad}, it follows that
\begin{corollary}
Let $W$ be a set consisting of the following
polynomials:
\begin{align*}
&[x_ix_j]_2-[x_jx_i]_1,\ \ [x_ix_j]_1-[x_jx_i]_2, \ \ [x_ix_i]_2-[x_ix_i]_1, &  (&i,j\in I,\ i>j).
\end{align*}
Then
$W$ is the reduced Gr\"{o}bner-Shirshov basis for the free abelian  dialgebra $Di\langle X\mid T\rangle$.
\end{corollary}

From Theorem \ref{tndsg} and Proposition \ref{pfad}, it follows that

\begin{corollary}(\cite[Theorem 1]{Zhu15})
Let
\begin{align*}
FAd(X):=\{\lfloor x_{i_1}\cdots x_{i_{t-1}}x_{i_t}x_{i_{t+1}}\cdots x_{i_l}\rfloor_t\mid &
x_{i_1}\leq\cdots\leq x_{i_{t-1}}<x_{i_t}\leq x_{i_{t+1}}\leq\cdots \leq x_{i_l}, \\
&x_{i_j}\in X,\ 1\leq j\leq l,\ \  l,t\in \mathbb{Z}^+, \ t\leq l
\}.
\end{align*}
Then
$( FAd(X), \ \vdash, \dashv \ )$ is the free abelian disemigroup generated by $X$,
where the operations $\vdash$ and $\dashv$ are as follows: for any $\lfloor u\rfloor_t, \lfloor v\rfloor_p\in FAd(X)$,
\begin{eqnarray*}
\lfloor u\rfloor_t\vdash\lfloor v\rfloor_p=\lfloor uv\rfloor_{_{\tau_{_{_{\lfloor v\rfloor \lfloor u\rfloor}}}(p)}}, \ \
\lfloor u\rfloor_t\dashv\lfloor v\rfloor_p=\lfloor uv\rfloor_{_{\tau_{_{_{\lfloor u\rfloor \lfloor v\rfloor}}}(t)}}.
\end{eqnarray*}
\end{corollary}

\subsection{Normal forms of free left (right) commutative disemigroups}

\begin{definition}
\emph{A disemigroup} $(D,\vdash,\dashv)$ \emph{is} left (right)  commutative
\emph{if} $a\dashv b\dashv c=b\dashv a\dashv c, \ a\vdash b\vdash c=b\vdash a\vdash c$ \
$(a\vdash b\vdash c=a\vdash c\vdash b, \ a\dashv b\dashv c=a\dashv c\dashv b)$ \emph{for all} $a,b,c\in D$.

\end{definition}

Let $X$ be an arbitrary set and $T$ the subset of $[X^+]_\omega\times[X^+]_\omega$,
where $T$ consists of the following:
$$
([u]_m\vdash [v]_n\vdash [w]_l,[v]_n\vdash [u]_m\vdash [w]_l), \ \  ([u]_m\dashv [v]_n\dashv [w]_l,[v]_n\dashv [u]_m\dashv [w]_l),
$$
where $[u]_m,[v]_n,[w]_l\in [X^+]_\omega$. Then $Disgp\langle X\mid T\rangle$  is the free left commutative disemigroup generated by $X$.

\begin{proposition}\label{pflcd}
Let $X=\{x_i\mid i\in I\}$ be a well-ordered set and $T$ the subset of $Di\langle X\rangle$
consisting of the following polynomials:
\begin{eqnarray*}
[u]_m\vdash [v]_n\vdash [w]_l-[v]_n\vdash [u]_m\vdash [w]_l,\ \
[u]_m\dashv [v]_n\dashv [w]_l-[v]_n\dashv [u]_m\dashv [w]_l,
\end{eqnarray*}
where $[u]_m,[v]_n,[w]_l\in [X^+]_\omega$. Then
\begin{enumerate}
\item[(i)] \ $Di\langle X\mid T\rangle=Di\langle X\mid S\rangle$, where $S$ consists of the following
polynomials:
\begin{align*}
&[uxv]_{|u|+1}-[\lfloor u\rfloor xv]_{|u|+1}, &  (&u,v\in X^*,x\in X, \ |u|\geq 2,|v|\leq 1),\\
&[uxvy]_{|u|+1}-[\lfloor uxv\rfloor y]_1, &  (&u,v\in X^*, x,y\in X,\ |v|\geq1).
\end{align*}
\item[(ii)] \ $S$ is a Gr\"{o}bner-Shirshov basis in $Di\langle X\rangle$.
\item[(iii)] \ The set
\begin{align*}
& \{[x_{i_1}\dots x_{i_n}]_{1}\mid x_{i_1}\leq\cdots\leq x_{i_{n-1}}, \ x_{i_l}\in X, \ 1\leq l \leq n, \ n\in \mathbb{Z}^+ \}\\
& \cup\{[x_{j_1}\dots x_{j_m}u]_{m} \mid x_{j_1}\leq\cdots\leq x_{j_{m-1}},\ x_{j_k}\in X,\ 1\leq k \leq m, \
m\in \mathbb{Z}^+, u\in X^*, |u|\leq 1 \}
\end{align*}
is a $\mathbf{k}$-basis of the dialgebra $Di\langle X\mid T\rangle$ and
normal forms of elements of the free left commutative disemigroup $Disgp\langle X\mid T\rangle$.
\end{enumerate}
\end{proposition}
\noindent{\bf Proof.}
$(i)$ We only need to prove that the polynomials in $S$ are trivial modulo $T$ and the polynomials in $T$ are trivial modulo $S$.
\begin{align*}
[uxv]_{|u|+1}-[\lfloor u\rfloor xv]_{|u|+1}&=[u]_{|u|}\vdash [xv]_1-\lfloor u\rfloor_{|u|}\vdash [xv]_1 \\
& \equiv  \lfloor u\rfloor_{|u|}\vdash [xv]_1-\lfloor u\rfloor_{|u|}\vdash [xv]_1\equiv 0 \ \ \mathrm{mod} (T),\\
 [uxvy]_{|u|+1}-[\lfloor uxv\rfloor y]_1&=[ux]_{|u|+1}\dashv [v]_1\dashv y-\lfloor uxv\rfloor_1\dashv y \\
& \equiv  [v]_1\dashv [ux]_{|u|+1} \dashv y -\lfloor uxv\rfloor_1\dashv y\\
&\equiv [vux]_1\dashv y-\lfloor uxv\rfloor_1\dashv y \\
& \equiv  \lfloor vux\rfloor_1\dashv y-\lfloor uxv\rfloor_1\dashv y \equiv 0 \ \ \mathrm{mod} (T).
\end{align*}
Suppose that $[w]_l=[w_1]_{|w_1|}\vdash [xw_2]_1=[w'y]_l$
where $w_1,w_2,w'\in X^*, x,y\in X$.
\begin{eqnarray*}
 \ \ \ \ \ \ \  \ \ \ \ \ \ \ \ \ [u]_m\vdash [v]_n\vdash [w]_l-[v]_n\vdash [u]_m\vdash [w]_l
 =[uvw_1 x w_2]_{|uvw_1|+1}-[vuw_1 x w_2]_{|vuw_1|+1}
\end{eqnarray*}
\begin{displaymath}  \ \ \ \ \ \ \ \ \ \ \ \ \ \ \ \ \ \  \equiv
\begin{cases}
[\lfloor uvw_1 \rfloor x w_2]_{|uvw_1|+1}-[\lfloor vuw_1 \rfloor x w_2]_{|vuw_1|+1}\equiv 0 \ \ \mathrm{mod} (S)
 \ \ \text{if $|w_2|\leq 1$,}\\
[\lfloor uvw_1 x w'_2\rfloor y]_1- [\lfloor vuw_1 x w'_2\rfloor y]_1 \equiv 0  \  \mathrm{mod} (S) \
\text{if $|w_2|> 1, w_2=w'_2y$.}
\end{cases}
\end{displaymath}
\begin{align*}
[u]_m\dashv [v]_n\dashv [w]_l-[v]_n\dashv [u]_m\dashv [w]_l=&[uvw'y]_{m}-[vuw'y]_{n} \\
\equiv &[\lfloor uvw'\rfloor y]_1 -[\lfloor vuw'\rfloor y]_1 \equiv 0  \  \mathrm{mod} (S).
\end{align*}

$(ii)$ It is easy to check that all possible compositions of left (right) multiplication in $S$ are equal to zero.
For any composition of $(f,g)_{[w]_m}$ in $S$, note that $-r\!_{_f},-r\!_{_g}\in [X^+]_\omega$, $|w|\geq 4$,
$[w]_m=\overline{[afb]_m}=\overline{[cgd]_m}$ and
$\lfloor w\rfloor=\lfloor a\widetilde{r\!_{_f}}b\rfloor=\lfloor c\widetilde{r\!_{_g}}d\rfloor$,
where $a,b,c,d\in X^*$. It follows that
\begin{eqnarray*}
(f,g)_{[w]_m}=[afb]_m-[cgd]_m=-[a\widetilde{r\!_{_f}}b]_{m_1}+[c\widetilde{r\!_{_g}}d]_{m_2}.
\end{eqnarray*}
If $|w|-m\leq 1$, then $m_1=m_2=m$ and $|a\widetilde{r\!_{_f}}b|-m_1\leq 1$, $|c\widetilde{r\!_{_g}}d|-m_2\leq 1$.
\begin{eqnarray*}
-[a\widetilde{r\!_{_f}}b]_{m_1}+[c\widetilde{r\!_{_g}}d]_{m_2}\equiv -[\lfloor u\rfloor xv]_{|u|+1}+[\lfloor u\rfloor xv]_{|u|+1}
\equiv 0 \ \ \mathrm{mod} (S),
\end{eqnarray*}
where $u,v\in X^*$, $x\in X$, $|u|\geq 2, |v|\leq1$. \\
If $|w|-m>1$, then $m_1\leq m, m_2\leq m$ and $|a\widetilde{r\!_{_f}}b|-m_1> 1$, $|c\widetilde{r\!_{_g}}d|-m_2> 1$.
\begin{eqnarray*}
-[a\widetilde{r\!_{_f}}b]_{m_1}+[c\widetilde{r\!_{_g}}d]_{m_2}\equiv -[\lfloor uxv\rfloor y]_1+[\lfloor uxv\rfloor y]_1
\equiv 0 \ \  \mathrm{mod} (S),
\end{eqnarray*}
where $u,v\in X^*$, $x,y\in X$, $|v|\geq 1$.

Then all the compositions in $S$ are trivial. We have proved $(ii)$.

$(iii)$ This part follows from  Theorems \ref{cd} and \ref{tndsg}.
\ \ $\square$

From Theorem \ref{cd},  Lemma \ref{lirr} and Proposition \ref{pflcd}, it follows that
\begin{corollary}
Let $W$ be a set consisting of the following
polynomials:
\begin{align*}
&[x_ix_jx_t]_3-[x_jx_ix_t]_3,,\ \ [x_ix_jx_t]_1-[x_jx_ix_t]_1, & (&i,j,t\in I,\ i>j),\\
&[x_lx_ix_jx_t]_2-[\lfloor x_lx_ix_j\rfloor x_t]_1,  & (&l,i,j,t\in I,\ i\leq j).
\end{align*}
Then
$W$ is the reduced Gr\"{o}bner-Shirshov basis for the free left commutative dialgebra $Di\langle X\mid T\rangle$.
\end{corollary}

Analysis similar to that in the proof of Proposition \ref{pflcd} shows the following proposition.
\begin{proposition}\label{pfrcd}
Let $X=\{x_i\mid i\in I\}$ be a well-ordered set and $T'$ the subset of $Di\langle X\rangle$
consisting of the following polynomials:
\begin{eqnarray*}
[w]_l\dashv [v]_n\dashv [u]_m-[w]_l\dashv [u]_m\dashv [v]_n,\ \
[w]_l\vdash [v]_n\vdash [u]_m-[w]_l\vdash [u]_m\vdash [v]_n,
\end{eqnarray*}
where $[u]_m,[v]_n,[w]_l\in [X^+]_\omega$. Then
\begin{enumerate}
\item[(i)] \ $Di\langle X\mid T'\rangle=Di\langle X\mid S'\rangle$, where $S'$ consists of the following
polynomials:
\begin{align*}
&[vxu]_{|v|+1}-[vx\lfloor u\rfloor]_{|v|+1}, &  (&u,v\in X^*,x\in X, \ |u|\geq 2,|v|\leq 1),\\
&[yvxu]_{|v|+2}-[y\lfloor vxu\rfloor]_{3}, &  (&u,v\in X^*, x,y\in X,\ ,|v|\geq1).
\end{align*}
\item[(ii)] \ $S$ is a Gr\"{o}bner-Shirshov basis in $Di\langle X\rangle$.
\item[(iii)] \ The set
\begin{align*}
&\{[ux_{i_1}\dots x_{i_n}]_{|u|+1}\mid x_{i_2}\leq\cdots\leq x_{i_{n}},\ x_{i_l}\in X, \  1\leq l\leq n, \ n\in \mathbb{Z}^+, \
 u\in X^*, |u|\leq 1  \}\\
& \cup\{[x_{j_1}\dots x_{j_{m}}]_{3} \mid x_{j_2}\leq\cdots\leq x_{j_{m}},\ x_{j_k}\in X,\ 1\leq k\leq m,\
m\in \mathbb{Z}^+, \ m\geq 3 \}
\end{align*}
is a $\mathbf{k}$-basis of the dialgebra $Di\langle X\mid T'\rangle$ and
normal forms of elements of the free right commutative disemigroup $Disgp\langle X\mid T'\rangle$.
\end{enumerate}
\end{proposition}

\begin{corollary}
Let $W'$ be a set consisting of the following
polynomials:
\begin{align*}
&[x_tx_ix_j]_1-[x_tx_jx_i]_1,\ \ [x_tx_ix_j]_3-[x_tx_jx_i]_3, & (&t,i,j\in I,\ i>j),\\
&[x_tx_ix_jx_l]_3-[x_t\lfloor x_ix_jx_l\rfloor]_3,  & (&t,i,j,l\in I,\ i\leq j).
\end{align*}
Then
$W'$ is the reduced Gr\"{o}bner-Shirshov basis for the free right commutative dialgebra $Di\langle X\mid T'\rangle$.
\end{corollary}


\begin{thebibliography}{10}

\bibitem{AL}
W.~W. Adams, P.~Loustaunau, An introduction to {G}r\"obner bases,
  Graduate Studies in Mathematics, vol. 3,  American Mathematical
  Society, 1994.

\bibitem{be78}
G.~M. Bergman, The diamond lemma for ring theory, Adv. in Math. 29(2) (1978)
  178--218.

\bibitem{bo76}
L.~A. Bokut, Imbeddings into simple associative algebras, Algebra i Logika
  15(2) (1976) 117--142, 245.

\bibitem{BC}
L.~A. Bokut, Y.~Chen, Gr\"obner-{S}hirshov bases: Some new results, in: Proceedings of the Second International Congress
in Algebra and Combinatorics, World Scientific, 2008, pp. 35--56.

\bibitem{BokutChenBook}
L.~A. Bokut, Y.~Chen, Gr\"obner-Shirshov Bases and Shirshov Algorithm,
  Educational tutorial lecture notes, Novosibirsk State University, 2014.

\bibitem{Di}
L.~A. Bokut, Y.~Chen, C.~Liu, Gr\"obner-{S}hirshov bases for dialgebras,
  Internat. J. Algebra Comput. 20(3) (2010) 391--415.

\bibitem{BC14}
L.~A. Bokut, Y.~Chen, Gr\"obner-{S}hirshov bases and their calculation, Bull.
  Math. Sci. 4(3) (2014) 325--395.

\bibitem{BoFKK00}
L.~A. Bokut, Y.~Fong, V.-F. Ke, P.~S. Kolesnikov, Gr\"obner and
  {G}r\"obner-{S}hirshov bases in algebra, and conformal algebras, Fundam.
  Prikl. Mat. 6(3) (2000) 669--706.

\bibitem{BK03}
L.~A. Bokut, P.~S. Kolesnikov, Gr\"obner-{S}hirshov bases: from inception to
  the present time, Zap. Nauchn. Sem. S.-Peterburg. Otdel. Mat. Inst. Steklov.
  (POMI) 272(Vopr. Teor. Predst. Algebr i Grupp. 7) (2000) 26--67, 345.

\bibitem{BK05a}
L.~A. Bokut, P.~S. Kolesnikov, Gr\"obner-{S}hirshov bases, conformal algebras,
  and pseudo-algebras, Sovrem. Mat. Prilozh. (13, Algebra) (2004) 92--130.

\bibitem{BKu94}
L.~A. Bokut, G.~P. Kukin, Algorithmic and combinatorial algebra,
  Mathematics and its Applications, vol. 255, Kluwer Academic
  Publishers Group, Dordrecht, 1994.


\bibitem{bu70}
B.~Buchberger, Ein algorithmisches {K}riterium f\"ur die {L}\"osbarkeit eines
  algebraischen {G}leichungssystems, Aequationes Math. 4 (1970) 374--383.

\bibitem{BuCL}
B.~Buchberger, G.~Collins, R.~Loos, R.~Albrecht, Computer algebra, symbolic and
  algebraic computation, Computing Supplementum, vol. 4,
  Springer-Verlag, New York, 1982.

\bibitem{BuW}
B.~Buchberger, F.~Winkler, Gr\"obner bases and applications, London
  Mathematical Society Lecture Note Series, vol. 251, Cambridge University
  Press, Cambridge, 1998.

\bibitem{CLO}
D.~A. Cox, J.~Little, D.~O'Shea, Ideals, varieties, and algorithms: An
  introduction to computational algebraic geometry and commutative algebra,
  Undergraduate Texts in Mathematics, Springer, Cham, 4th edition, 2015.

\bibitem{Ei}
D.~Eisenbud, Commutative algebra: With a view toward algebraic geometry, Graduate Texts in Mathematics,
vol. 150, Springer-Verlag, New York, 1995.

\bibitem{Fe}
R.~Felipe, An analogue to functional analysis in dialgebras, Int. Math. Forum
  2(21-24) (2007) 1069--1091.

\bibitem{Hi64}
H.~Hironaka, Resolution of singulatities of an algebtaic variety over a field
  if charac-teristic zero, {I}, {II}, Math. Ann. 79 (1964) 109--203, 205--326.

\bibitem{Ko}
P.~S. Kolesnikov, Varieties of dialgebras, and conformal algebras, Sibirsk.
  Mat. Zh. 49(2) (2008) 322--339.


\bibitem{Lo95}
J.-L. Loday, Alg\`ebres ayant deux op\'erations associatives (dig\`ebres), C.
  R. Acad. Sci. Paris S\'er. I Math. 321(2) (1995) 141--146.

\bibitem{Lo99}
J.-L. Loday, A.~Frabetti, F.~Chapoton, F.~Goichot, Dialgebras and related
  operads, Lecture Notes in Mathematics, vol. 1763, Springer-Verlag,
  Berlin, 2001.


\bibitem{Nm}
M.~H.~A. Newman, On theories with a combinatorial definition of
  ``equivalence.'', Ann. of Math. (2) 43 (1942) 223--243.

\bibitem{Po}
A.~P. Pozhidaev, 0-dialgebras with bar-unity and nonassociative {R}ota-{B}axter
  algebras, Sibirsk. Mat. Zh. 50(6) (2009) 1356--1369.

\bibitem{Sh62b}
A.~I. Shirshov, Some algorithmic problems for {$\varepsilon $}-algebras,
  Sibirsk. Mat. \v Z. 3 (1962) 132--137.

\bibitem{Shir3}
A.~I. Shirshov, Selected works of {A}. {I}. {S}hirshov, Translated from the
  Russian by Murray Bremner and Mikhail V. Kotchetov, Edited by Leonid A.
  Bokut, Victor Latyshev, Ivan Shestakov and Efim Zelmanov, Contemporary
  Mathematicians, Birkh\"auser Verlag, Basel, 2009.

\bibitem{Zhuchok}
A.~V. Zhuchok, Free commutative dimonoids, Algebra Discrete Math. 9(1) (2010)
  109--119.

\bibitem{Zhuchok17}
A.~V. Zhuchok, Structure of relatively free dimonoids, Comm. Algebra. 45(4) (2017)  1639--1656.



\bibitem{Zhu15}
Y.~V. Zhuchok, Free abelian dimonoids, Algebra Discrete Math. 20(2) (2015)
  330--342.

\end{thebibliography}
\end{document}